\RequirePackage{fix-cm}
\documentclass[smallextended]{svjour3}       
\smartqed  
\usepackage{graphicx}
\usepackage{bm}
\usepackage{amsmath,amssymb, amsfonts,mathtools}
\usepackage{epstopdf}
\usepackage{epsfig}
\usepackage{subcaption}
\usepackage{mathrsfs}
\usepackage{enumitem}

\usepackage{mathtools}
\mathtoolsset{showonlyrefs}
\begin{document}

\title{A structure-preserving scheme for computing effective diffusivity and anomalous diffusion phenomena of random flows
}
%
\titlerunning{A structure-preserving scheme for computing effective diffusivity of random flows}        

\author{ Tan Zhang \and Zhongjian Wang \and Jack Xin \and Zhiwen Zhang 
}

\authorrunning{Tan Zhang, Zhongjian Wang, Jack Xin and Zhiwen Zhang} 
\institute{Tan Zhang \at
              Department of Mathematics, The University of Hong Kong, Pokfulam, Hong Kong. \\
              \email{thta@connect.hku.hk}           
           \and
           Zhongjian Wang \at
           Division of Mathematical Sciences, School of Physical and Mathematical Sciences, Nanyang  Technological University, Singapore 637371.\\
           \email{zhongjian.wang@ntu.edu.sg}
           \and
           Jack Xin \at
           Department of Mathematics, University of California at Irvine, Irvine, CA 92697, USA. \\
           \email{jxin@math.uci.edu}
           \and
           Zhiwen Zhang \at
           Department of Mathematics, The University of Hong Kong, Pokfulam, Hong Kong. \\
           \email{zhangzw@hku.hk} 
}
\date{Received: date / Accepted: date}

\maketitle
\begin{abstract}
\noindent This paper aims to investigate the diffusion behavior of particles moving in stochastic flows under a structure-preserving scheme. We compute the effective diffusivity for normal diffusive random flows and establish the power law between spatial and temporal variables for cases with anomalous diffusion phenomena.
From a Lagrangian approach, we separate the corresponding stochastic differential equations (SDEs) into sub-problems and construct a one-step structure-preserving method to solve them. 
Then by modified equation systems, the convergence analysis in calculating the effective diffusivity is provided and compared between the structure-preserving scheme and the Euler-Maruyama scheme.
Also, we provide the error estimate for the structure-preserving scheme in calculating the power law for a series of super-diffusive random flows. Finally, we calculate the effective diffusivity and anomalous diffusion phenomena for a series of 2D and 3D random fields. 
\end{abstract}

\section{Introduction}
\label{sec:intro}
\noindent The study of diffusion enhancement phenomena in fluid advection has been ongoing since the groundbreaking research conducted by Sir G. Taylor\cite{Taylor1922}. A key challenge lies in accurately characterizing and quantifying the effective diffusion at large scales in complex fluid flows, as this holds significant theoretical and practical significance. The applications of this research span various scientific and technological domains, including atmospheric science, oceanography, chemical engineering, and combustion. 

In this work, we aim to compute the effective diffusivity of the passive tracer models and investigate anomalous diffusion
phenomena in chaotic and random flows. For the passive tracer models in random flows, it can be described by the following stochastic differential equations (SDEs).
\begin{equation}
    d\bm x(t) = \bm v(\bm x, t, \omega) dt + \sigma d\bm w(t),\ \bm x \in \mathbb{R}^d,\ t \geq 0,
    \label{PasTraSDEforVx}
\end{equation}
where $\bm x(t) = (x_1(t),...,x_d(t))^T \in \mathbb{R}^d$ is the spatial position of the passive tracer particle, $\bm v(\bm x, t, \omega)$ is the incompressible velocity field, $\sigma > 0$ is the molecular diffusivity, and $\bm w(t)$ is the standard $d$-dimensional Brownian motion. The velocity field $\bm v(\bm x, t, \omega)$ is a mean zero, stationary, ergodic random vector field over a certain probability space. Here $\omega$ is an element of the probability space describing all possible environments which is exclusive of the Brownian motion, $\bm w(t)$. We assume that the realizations of $\bm v(\bm x, t, \omega)$ are almost surely locally Lipschitz in $\bm x$. This is to guarantee the existence of the solution of equation \eqref{PasTraSDEforVx}. In practice, for a given energy spectrum density, we follow the generation method in \cite{Kraichnan1970}\cite{Majda1999} to model the random field $\bm v(\bm x, t, \omega)$.

One of the primary challenges associated with the SDEs \eqref{PasTraSDEforVx} involves investigating the behavior of the particles $\bm x(t)$ in \eqref{PasTraSDEforVx} over an extended time and on a large scale. For the rescaled process, $\epsilon \bm x(t / \epsilon^2)$, it is crucial to verify whether it will converge in distribution to a new Brownian motion $\Tilde{\bm w}(t)$ with a covariance matrix $\mathbf{D}^E$. If so, then $\mathbf{D}^E$ is denoted as the effective diffusivity matrix. This problem can be viewed as the homogenization of the passive tracer model and is highly nontrivial.  Extensive research has been conducted to compute the effective diffusivity matrix $\mathbf{D}^E$ under various conditions on the flows. In situations where there are spatial-temporal periodic fields or random fields with short-range correlations, the homogenization theory can be applied to calculate the effective diffusivity matrix $\mathbf{D}^E$ \cite{Panpanicolaou2011},\cite{Garnier1997},\cite{Jikov1994},\cite{Pavliotis2008}. This enables us to determine the behavior of the diffusion process in such systems. 

In the case of random flows characterized by long-range correlations, studying the behavior of particle motion over an extended time and on a large scale poses significant challenges due to its intricate and general complexity. The presence of various forms of anomalous diffusion, including super-diffusion and sub-diffusion, further adds to the difficulty of the analysis.

Multiple theoretical studies have been conducted on the homogenization of time-dependent random flows \cite{Carmona1997},\cite{Fannjiang1999},\cite{Fannjiang1996},\cite{Komorowski2001},\cite{Landim1998}. One notable example of these results is the work by the authors in \cite{Carmona1997}, where they successfully demonstrated the existence of the effective diffusivity for the 2D incompressible time-dependent Gaussian field. In \cite{Fannjiang1999},\cite{Fannjiang1996}, the authors demonstrated certain necessary conditions in turbulent flows that should be satisfied for the behavior of advection-diffusion over an extended time to exhibit normal diffusion behavior. For random flows that exhibit normal diffusion behavior, we have proposed a structure-preserving method to compute the effective diffusivity in the Lagrangian approach \cite{Junlong2021}. In our previous studies, the focus has primarily been on random flows that exhibit diffusive behavior, rather than on random fields that display abnormal diffusion. Moreover, the numerical analysis for calculating the effective diffusivity is unexplored when the diffusion term's coefficient $\sigma$ vanishes.
This motivates us to provide the following contributions in this paper: 
\begin{itemize}[label=$\bullet$]
    \item Provide the convergence analysis of our structure-preserving method in computing effective diffusivity when the diffusion coefficient $\sigma$ vanishes;
    \item Compare the error estimate of our structure-preserving method with the traditional Euler-Maruyama scheme;
    \item Investigate the anomalous diffusion phenomena and verify the accuracy of our structure-preserving method in a series of random flows.
\end{itemize}

The whole procedure is as follows. We first apply an operator splitting method to separate the SDE system \eqref{PasTraSDEforVx} into two sub-problems. In our structure-preserving method, we adopt an implicit structure-preserving scheme to solve the advection term and adopt the Euler-Maruyama scheme to solve the diffusion term. 
From the perspective of the modified equation, we obtain the continuous modified equation systems \eqref{ModifiedEquofEul} which correspond to our structure-preserving method and Euler-Maruyama scheme respectively. Then by applying Ito's formula, we can analyze and obtain the error estimate of two modified equation systems in computing effective diffusivity (Theorem \ref{Thm::ErrEstforBotnSch}). We have also analyzed that the structure-preserving method's modified system can maintain the diffusion behavior of the original system better when $\sigma$ is relatively small, regardless of whether it is sub-diffusive (Theorem \ref{Coro::subdiffforSP}) or normal diffusive (Corollary \ref{Coro::diffSP}). While the Euler-Maruyama scheme does not hold this property. Meanwhile, for a group of random velocity fields that exhibit super-diffusion as mentioned in \cite{Fannjiang1998}, we prove that our structure-preserving method can maintain its original super-diffusion behavior (Theorem \ref{ThmforSuperdiff}). 
Finally, we present numerical results to verify the theorems, demonstrate the convection-enhanced diffusion phenomenon, and compute the effective diffusivity and anomalous diffusion phenomena in 2D and 3D random flows.

The rest of the paper is organized as follows. In Section 2, we introduce the definition of effective diffusivity and briefly review some existing results for diffusion in random flows. We also provide the sharp condition for random fields to exhibit normal diffusive phenomena. In Section 3, we develop the implicit structure-preserving method for computing the effective diffusivity of passive tracer models in random flows. We also make a comparison of the volume-preserving property between the structure-preserving method and the Euler-Maruyama method. Moreover, we introduce the modified equation systems of two numerical schemes. In Section 4, we provide the convergence analysis for our structure-preserving method and Euler-Maruyama method based on the modified equation system. Error analysis of the performance of the structure-preserving scheme under random flows with super-diffusive behavior is also carried out. In Section 5, we conduct numerical experiments and present results to verify the error analyses we proved, demonstrate the convection-enhanced diffusion phenomenon, and compute the effective diffusivity and anomalous diffusion phenomena in 2D and 3D random flows.

\section{Preliminaries}
To ensure the self-contained nature of this paper, we provide a concise overview of the current research in the area of advection-enhanced diffusion in stochastic turbulent flows and the definition of effective diffusivity of the associate passive tracer models \eqref{PasTraSDEforVx}. Some results for the general cases we put in Appendix \ref{AppendixforSomeFormandSGC}, and in the main body we use these results in a more specific system (Subsection \ref{subsec::CorProbandDefogE}). We will also go through the criterion to determine whether a stochastic turbulent flow will perform normal diffusive behavior or not (Subsection \ref{subsec::Sharpcondi}).


\subsection{Continuous corrector problem and the definition of effective diffusivity}
\label{subsec::CorProbandDefogE}
 To study the effective diffusivity of the stochastic flows connected to \eqref{PasTraSDEforVx}, we first present a concise overview of the current findings regarding convection-enhanced diffusion in random flows.

  Let $(\mathcal{X}, \mathcal{H}, P_0)$ be a probability space. Here $\mathcal{H}$ is a $\sigma$-algebra over $\mathcal{X}$, $P_0$ is a probability measure. Let $\tau_{\bm x}$, $\bm x \in \mathbb{R}^d$, be an almost surely continuous, jointly measurable group of measure-preserving transformation on $\mathcal{X}$. The properties of this system are illustrated in Appendix \ref{AppendixforSomeFormandSGC} and here we follow its defined notations. Then let $\Omega$ be the space of $\mathcal{X}$-valued continuous function $C([0,\infty); \mathcal{X})$ and let $\mathcal{B}$ be its Borel $\sigma$-algebra. Then we denote $\omega(\cdot + t), t \geq 0$ as the standard shift operator on the path space $(\Omega, \mathcal{B})$.

   From a Lagrangian particle approach, we can define the following environment process  $\eta: [0, \infty)\times \Omega \to \mathcal{X}$ at time $t$ as 
\begin{equation}
    \eta(t) = \tau_{\bm x_t^{\omega}} \omega(t),\ \text{where}\ \eta(0) = \omega(0).
    \label{ContinTimeCase} 
\end{equation}
Here $\bm x_t^{\omega} \in \mathbb{R}^d$ is the position of the particle at time $t$ with the realization of the environment $\omega$.  The random field $\bm v(\bm x, t, \omega)$ is stationary provided $\tau_{\bm x}$ preserves the probability measure $P_0$, for every $\bm x \in \mathbb{R}^d$. And the ergodicity of $\bm v(\bm x, t, \omega)$ is also provided by the transformation $\tau_{\bm x}, \bm x \in \mathbb{R}^d$ is ergodic. Furthermore, inconvenience to the subsequent analysis of the effective diffusivity, we require that the incompressible random velocity field $\bm v(\bm x, t, \omega)$ satisfies the following assumptions.
 \begin{assumption}
     \label{Assforvpro}
      The incompressible random field $\bm v(\bm x, t, \omega)$ is jointly continuous in $(\bm x, t)$, bounded and with finite second moments.
 \end{assumption}

 \begin{assumption}
\label{Assforvregular}
     $\bm v(\bm x, t, \omega)$ is locally Lipschitz in the spatial variable and has certain regularity, i.e. $ \bm v \in (C_b^m(\mathcal{X}))^d$ for some $m \geq 1$.
\end{assumption}

 By denoting $\textbf{M}$ as the expectation corresponding to the Brownian motion term in the SDEs, the environment process will generate a semigroup of transformation
\begin{equation}
     \bm S^t f(\Tilde{\bm x}) = \textbf{ME}_{\Tilde{\bm x}} f(\eta(t)), t\geq 0,\ f\in L^{\infty}(\mathbb{R}^d),
\end{equation}

For clarity, we denote $\mathbb{E}$ as the total expectation with respect to all randomness, i.e., $\mathbb{E} = \mathbf{M E}$, in the remaining part of this paper.
From Proposition 3 and Proposition 6 in \cite{Fannjiang1999}, we know that $\bm S^t$ is a strongly continuous, measure-preserving, Markov semigroup of contraction on $L^2(\mathbb{R}^d)$.  
Then we let $\mathcal{L}_0$ denote the generator of the semigroup $\bm S^t, t\geq 0$, i.e. 

\begin{equation}
    \mathcal{L}_0 f = \frac{\sigma^2}{2} \Delta f + \bm v \cdot \nabla f.
    \label{defofL0}
\end{equation}
$\mathcal{L}_0$ is a concretization of generator $L$ which is mentioned in Appendix \ref{AppendixforSomeFormandSGC}. So $\bm S^t$ also holds the spectral gap condition \eqref{SpectralGapProp}, i.e. 
\begin{equation}
    ||\bm S^t f||_{L^2(\mathbb{R}^d)}\leq \exp(-2c_1 t)||f||_{L^2(\mathbb{R}^d)},\ c_1 > 0,\ f\in L_0^2(\mathbb{R}^d).
    \label{SpectralGapPropforL}
\end{equation}

The proof of \eqref{SpectralGapPropforL} can be seen from Proposition 6 in \cite{Fannjiang1999}. Then we can define $\bm G = \int_0^{\infty} \bm S^t \bm v dt$ and it satisfies the following continuous corrector problem, 
\begin{equation}
    \mathcal{L}_0 \bm G = - \bm v.
    \label{CorrectorProbRanContin}
\end{equation}

By solving the continuous corrector problem \eqref{CorrectorProbRanContin}, we can define the effective diffusivity.
\begin{proposition}
    [\cite{Fannjiang1999}, Lemma 1]Let $\bm x(t)$ be the solution to the passive tracer model \eqref{PasTraSDEforVx} and $\bm x_{\epsilon} (t) \equiv \epsilon \bm x(t / \epsilon^2)$. For any unit vector $\bm e \in \mathbb{R}^d$, let $\bm G_{\bm e}$ denote the projection of the vector solution $\bm G$ in \eqref{CorrectorProbRanContin} along the direction $\bm e$. Then the law of the process of the projection of $\bm x_{\epsilon}(t)$ along direction $\bm e$, i.e. $\bm x_{\epsilon}(t) \cdot \bm e$ will weakly converge to a Brownian motion with a diffusion coefficient which is given by 
    \begin{equation}
        \bm e^T \mathbf{D}^E \bm e = D_0 + (-\mathcal{L}_0\bm G_{\bm e}, \bm G_{\bm e})_{L^2(\mathbb{R}^d)} ,
        \label{DefofEffDiff}
    \end{equation}
    where $\mathbf{D}^E$ is the effective diffusivity associated with \eqref{PasTraSDEforVx} and $D_0 = \sigma^2 / 2$. 
\end{proposition}

We can see that the correction to $D_0$ in \eqref{DefofEffDiff} is non-negative definite and for any unit vector $\bm e \in \mathbb{R}^d$, $\bm e^T \mathbf{D}^E \bm e \geq D_0$. This is called convection-enhanced diffusion. The effective diffusivity fundamentally captures the long-time, large-scale behaviors of passive tracer models. In practice, as $D_0$ approaches extremely small values, the solutions of \eqref{CorrectorProbRanContin} exhibit sharp gradients, requiring a significant number of Fourier modes to resolve. Also here we consider the flow with unbounded domain setting. These characteristics render the Eulerian approach computationally demanding. Alternatively, the effective diffusivity matrix can be computed using the Lagrangian approach, defined by,
\begin{equation}
    D^E_{ij} = \lim_{t\to \infty} \frac{\big\langle \big(x_i(t) - x_i(0)\big) \big(x_j(t) - x_j(0)\big)\big\rangle}{2t},\ 1 \leq i,j \leq d,
    \label{DELagranAppro}
\end{equation}
where $\langle \cdot\rangle$ is the average taken over a collection of particles. If $D_{ij}^E$ exists and is a nonzero value, then it illustrates that the transport of the particle under this random field exhibits a diffusive behavior on a long-time scale.

If $\bm v(\bm x, t, \omega)$ satisfies Assumption \ref{Assforvpro} and \ref{Assforvregular}, then the dispersion of the particles $\big\langle \big(x_i(t) - x_i(0)\big) \big(x_j(t) - x_j(0)\big)\big\rangle$ will grow linearly with respect to the time $t$. In our previous study \cite{Junlong2021}, we have given the convergence analysis of the effective diffusivity calculated by the structure-preserving scheme for random fields with normal diffusive performance. 

There are also instances where the dispersion of the transport of the particle does not exhibit linear growth with respect to the time $t$, but instead follows a power law of $\gamma$, i.e. $\mathbb{E}[x_i^2(t)] \sim t^{\gamma}$, where $1<\gamma <2$ and $0<\gamma < 1$ are corresponding to super-diffusion and sub-diffusion respectively. Such random fields with abnormal diffusion behavior are what we want to focus on in this paper, and we also want to compare the structure-preserving scheme with other methods of computing effective diffusivity. Due to our method being implemented by the Lagrangian particle approach, the computational cost is relatively small in the high-dimensional cases and when $D_0$ is extremely small. The major difficulty for the Lagrangian approach is to guarantee the accuracy during the evolution which is long enough to approach the diffusion time scale. This motivates us to develop a structure-preserving scheme that is accurate for long-time integration and compare it with the explicit Euler-Maruyama scheme in this problem.

\subsection{Sharp conditions for random flows to exhibit normal diffusive behavior}
\label{subsec::Sharpcondi}
In this subsection, we discuss the criterion of random fields to exhibit  diffusion behavior from their energy spectrum density. 
Notice that we have suppressed the dependence of the velocity on $\omega$ in $\bm v(\bm x, t)$ for notation simplicity here. In our following analysis and experiments, we will consider the particles moving in an isotropic velocity field $\bm v(\bm x, t)$ satisfies Assumption \ref{Assforvpro}. Thus, the random field $\bm v(\bm x, t)$ is fully specified by its covariance function $C(r, t) = Cov(\bm v(\bm x, t), \bm v(\bm 0,0))$ with $r = ||\bm x||_2$. We assume the covariance functions have the following forms
\begin{equation}
    C(r, t) = 2 D(t) \int_0^{\infty} E(k) \frac{\sin(kr)}{kr} dk,\ \text{for 3D flows},
    \label{CovFct3DS3}
\end{equation}
and
\begin{equation}
    C(r, t) = 2 D(t) \int_0^{\infty} E(k) J_0(kr) dk,\ \text{for 2D flows},
    \label{CovFct2DS3}
\end{equation}
where $D(t)$ is the time-correlation function, $E(k)$ is the energy spectrum function, and $J_0(\cdot)$ is the Bessel functions of the first kind. Given the covariance function \eqref{CovFct3DS3} or \eqref{CovFct2DS3}, we use the randomization method in \cite{Majda1999} to generate realizations of the random velocity field. Specifically, we generate the velocity field $\bm v(\bm x, t)$ with the form 
\begin{equation}
    \bm v(\bm x, t) = \frac{1}{\sqrt{N}} \sum_{n=1}^N \big(\bm u_n \cos(\bm k_n \cdot \bm x + \theta_n t) + \bm w_n \sin(\bm k_n \cdot \bm x + \theta_n t)\big).
    \label{vGene}
\end{equation}

For two-dimensional random flows, we have 
\begin{equation}
    \bm u_n = \xi_n \bm k_n^{\perp},\ \bm w_n = \zeta_n \bm k_n^{\perp},\ \bm k_n^{\perp} = (- k_n^2, k_n^1),\ n = 1,2,...,N,
    \label{uw2Dform} \notag
\end{equation}
where $\xi_n$ and $\zeta_n$ are independent 1D Gaussian variables. For three-dimensional random flows, we have 
\begin{equation}
    \bm u_n = \bm \xi_n \times \bm k_n,\ \bm w_n = \bm \zeta_n \times \bm k_n,\ \bm k_n = (k_n^1, k_n^2, k_n^3),\ n = 1,2,...,N,
    \label{uw3Dform} \notag
\end{equation}
where the samples $\bm \xi_n$ and $\bm \zeta_n$ are independent 3D Gaussian variables. Sample points for each component of $\bm k_n$ are independently chosen according to the energy spectrum function $E(k)$.

Here we assume that the velocity field is independent of the time $t$ in the following analysis, i.e. $D(t) = 1$.
For a 2D random velocity field with an isotropic energy spectrum $E(k)$, a sharp condition, which has been proved in \cite{Fannjiang1996}, to determine whether it will perform normal diffusion behavior is that 
\begin{equation}
    \int \frac{1}{k^2} \hat{C}(k_1, k_2) d k_1 dk_2 < \infty.
    \label{sharpCondi}
\end{equation}
Here $k = \sqrt{k_1^2 + k_2^2}$ and $\hat{C}(k_1, k_2)$ is the Fourier transform of the covariance function $C(x_1, x_2)$ defined in \eqref{CovFct2DS3}. Here $r = \sqrt{x_1^2 + x_2^2}$. Then together with the correlation function for 2D random flows \eqref{CovFct2DS3}, we can further obtain that
\begin{align}
    \hat{C}(k_1, k_2) &= 2 \int E(k') J_0(k' r) e^{-i(x_1k_1 + x_2 k_2)} dx_1 dx_2 d k' \notag \\
     &= 2 \int E(k') J_0(k' r) e^{-i k r \cos\theta} r d\theta dr d k' \notag \\
    & = 4\pi \int E(k') J_0(k' r) J_0(k r) r d r d k' \notag \\ 
    & = 4\pi \int E(k') \frac{1}{k'} \delta(k' - k) dk' = 4\pi E(k) \frac{1}{k}.
\end{align}

Then the sharp condition \eqref{sharpCondi} for 2D cases can be transformed into 
\begin{equation}
    \int \frac{1}{k^2} \hat{C}(k_1, k_2) dk_1 dk_2 = 2\pi \int \frac{1}{k^2} 4\pi E(k) \frac{1}{k} k d k = 8\pi^2 \int \frac{E(k)}{k^2} dk < \infty.
    \label{sharpCondi2D}
\end{equation}

A similar analysis can be done for 3D random flows. For a 3D random velocity field with energy spectrum $E(k)$, the sharp condition \eqref{sharpCondi2D} to determine whether it will perform normal diffusion behavior will now be transformed into 
\begin{align}
     \hat{C}(k_1,&k_2,k_3) = 2 \int_0^{\infty}\int_{\mathbb{R}^3} E(k') \frac{\sin(k'r)}{k'r} e^{-i(k_1x_1 + k_2x_2+ k_3x_3)} dx_1 dx_2 dx_3 dk'\notag \\
    &= 2 \int_0^{\infty} \int_0^{\infty} \int_0^{2\pi} \int_0^{\pi} E(k') \frac{\sin(k'r)}{k'r} e^{-i k r \sin \theta \cos \phi} r^2 \sin \theta d\theta d\phi dr dk' \notag \\
    &= 4\pi \int_0^{\infty} \int_0^{\infty} \int_0^{\pi} E(k') \frac{\sin(k'r)}{k'r} J_0(kr\sin\theta)r^2 \sin \theta d\theta dr dk' \notag \\
    &= 8\pi \int_0^{\infty} \int_0^{\infty} E(k') \frac{\sin(k'r)}{k'r} \frac{\sin(kr)}{kr} r^2 dr dk' \notag \\
    &= 8\pi \int_0^{\infty} \int_0^{\infty} E(k') \frac{\sin(k'r)}{k'} \frac{\sin(kr)}{k} dr dk' .
\end{align}

\begin{align}
    \int \frac{1}{k^2} \hat{C}(k_1, &k_2, k_3) dk_1 dk_2 dk_3 \notag \\
    & = 4\pi \int_0^{\infty} \frac{1}{k^2} (8\pi \int_0^{\infty} \int_0^{\infty} E(k') \frac{\sin(k'r)}{k'} \frac{\sin(kr)}{k} dr dk') k^2 dk \notag \\
    &= 32\pi^2 \int_0^{\infty} \int_0^{\infty} \int_0^{\infty} E(k') \frac{\sin(k'r)}{k'} \frac{\sin(kr)}{k} dr dk' dk\notag \\
    &= 16 \pi^3 \int_0^{\infty} \int_0^{\infty}  \frac{E(k')\sin(k'r)}{k'} dr dk' < \infty.
    \label{sharpCondi3D}
\end{align}

With the criterion of $2$-dimensional and $3$-dimensional normal diffusion obtained, we can make the following assumption about the random velocity fields for the subsequent study of diffusive and sub-diffusive cases.
\begin{assumption}
    The random velocity field $\bm v(\bm x, t)$ possesses an isotropic energy spectral density $E(k)$ which satisfies \eqref{sharpCondi2D} in 2D cases or \eqref{sharpCondi3D} in 3D cases.
    \label{AssforEk}
\end{assumption}
 
\section{Structure-preserving scheme for computing effective diffusivity}\label{sec:SPschemeandErrest}
In this section, we would like to introduce our structure-preserving scheme for $2$-dimensional and $d$-dimensional ($d \geq 3$) random flows (Subsection \ref{subsec::Constof2D} and \ref{subsec::constofnd}). Due to the difference in nature between $2$-dimensional and $d$-dimensional random flows, the conditions required for constructing the corresponding numerical schemes are not identical, so we discuss the construction of the numerical schemes for the $2$-dimensional and $d$-dimensional cases separately. Meanwhile, we verify the volume-preserving property of our method and make a comparison with the Euler-Maruyama scheme (Subsection \ref{Sec::CompSPandEul}). We also introduce the modified equation systems of two numerical schemes which will be used in the convergence analysis section (Subsection \ref{subsec::Modifiedequsys}).
\subsection{Construction of the scheme for $2$-dimensional random flows} 
\label{subsec::Constof2D}
 Here we first introduce the numerical scheme for computing the effective diffusivity of $2$-dimensional random flows. 
 For a given random realization of the environment, $\omega$, we first split the passive tracer model \eqref{PasTraSDEforVx} into two sub-problems based on the operator splitting method.
\begin{equation}
    d \bm x(t) = \bm v(\bm x(t), \omega) dt ,\ \text{and}\ \  d\bm x(t)  = \sigma d\bm w(t).
\end{equation}

By applying the Lie-Trotter splitting method and discretizing the sub-problems, we can develop our numerical schemes. Specifically, from time $t_n$ to $t_{n+1} = t_n + \Delta t$, where $\Delta t > 0, t_0 = 0$, assuming the numerical solution $\bm x_n$ for time $t_n$ is given, we discretize the subproblems and obtain
\begin{align}
    &\bm x_n^* = \bm x_n + \Delta t \bm v(\frac{\bm x_n^* + \bm x_n}{2}, \omega),
    \label{SubProb1Discret} \\
    &\bm x_{n+1} = \bm x_{n}^* + \sigma \sqrt{\Delta t} \bm w_n,
    \label{SubProb2Discret}
\end{align}
here $\bm w_n$ is a $d$-dimensional standard Gaussian random variable. By solving the numerical scheme \eqref{SubProb1Discret} \eqref{SubProb2Discret}, we obtain the numerical solution $\bm x_{n+1}$, which approximates the exact solution $\bm x(t)$ to the SDE system \eqref{PasTraSDEforVx} at time $t = t_{n+1}$. 

Here we can see that we apply the Euler midpoint rule to discrete the first sub-problem and obtain numerical scheme \eqref{SubProb1Discret}. One natural question arises: Is the numerical scheme \eqref{SubProb1Discret} volume-preserving? Then we can derive that
\begin{equation}
    \frac{\partial \bm x_n^*}{\partial \bm x_n} = \mathbf{I} + \Delta t \mathbf{D}\bm v(\frac{\bm x_n^* + \bm x_n}{2}, \omega) (\frac{1}{2} \frac{\partial \bm x_n^*}{\partial \bm x_n} + \frac{1}{2} \mathbf{I}),
    \label{PartialxsPatrialxn}
\end{equation}
where $\mathbf{I}$ is the identity matrix of size $d$, and $\mathbf{D}\bm v(\cdot)$ is the Jacobin matrix. Equation \eqref{PartialxsPatrialxn} implies that 
\begin{equation}
    \frac{\partial \bm x_n^*}{\partial \bm x_n} = \big (\mathbf{I} - \frac{\Delta t}{2} \mathbf{D}\bm v(\frac{\bm x_n^* + \bm x_n}{2}, \omega)\big )^{-1} \big (\mathbf{I} + \frac{\Delta t}{2} \mathbf{D}\bm v(\frac{\bm x_n^* + \bm x_n}{2}, \omega)\big).
\end{equation}

To ensure that the numerical scheme \eqref{SubProb1Discret} is volume-preserving, it is necessary to require the condition that the determinant of $\frac{\partial \bm x_n^*}{\partial \bm x_n}$ should equal to one, i.e. $\det(\frac{\partial \bm x_n^*}{\partial \bm x_n}) = 1$. Inspired by the approach in \cite{KangFeng1989}, here we denote $P(\lambda) = \det\big(\mathbf{D}\bm v(\frac{\bm x_n^* + \bm x_n}{2}, \omega) - \lambda \mathbf{I}\big)$ as the characteristic determinant of $\mathbf{D}\bm v(\frac{\bm x_n^* + \bm x_n}{2},\omega)$. Then we can derive that
\begin{equation}
    \det(\frac{\partial \bm x_n^*}{\partial \bm x_n}) = \frac{(\frac{\Delta t}{2})^d \det\big(\mathbf{D}\bm v(\frac{\bm x_n^* + \bm x_n}{2}, \omega) + \frac{2}{\Delta t} \mathbf{I}\big)} {(\frac{-\Delta t}{2})^d \det\big(\mathbf{D}\bm v(\frac{\bm x_n^* + \bm x_n}{2}, \omega) - \frac{2}{\Delta t} \mathbf{I}\big)} = (-1)^d \frac{P(-\frac{2}{\Delta t})}{P(\frac{2}{\Delta t})}.
    \label{CondiforVolPre}
\end{equation}

\begin{proposition}[\cite{KangFeng1989}, Theorem 24]
    When $d = 2$, the numerical scheme \eqref{SubProb1Discret} is volume-preserving if the velocity field $\bm v$ is divergence-free. 
\end{proposition}
\begin{proof}
    When $d=2$, equation \eqref{CondiforVolPre} will become into $\det(\frac{\partial \bm x_n^*}{\partial \bm x_n}) = P(-\frac{2}{\Delta t})/P(\frac{2}{\Delta t})$. Notice that 
    \begin{equation}
        P(\lambda) = \lambda^2 - \text{Tr}\big(\mathbf{D}\bm v(\frac{\bm x_n^* + \bm x_n}{2}, \omega)\big) \lambda + \det\big(\mathbf{D}\bm v(\frac{\bm x_n^* + \bm x_n}{2}, \omega)\big) = \lambda^2.
    \end{equation}

When $\bm v$ is divergence-free, we can further derive that $P(\lambda) = P(-\lambda)$. Then we can obtain that $\det(\frac{\partial \bm x_n^*}{\partial \bm x_n}) = 1$  which guarantees the scheme \eqref{SubProb1Discret} is volume-preserving in this case. 
\hfill $\square$
\end{proof}

\subsection{Construction of the scheme for $d$-dimensional separable random flows}
\label{subsec::constofnd}
For the commonly cases with $d \geq 3$, we now require the condition that $P(\lambda) = (-1)^d P(-\lambda)$. To satisfy this, it requires the velocity field $\bm v$ is divergence-free and $\det\big(\mathbf{D}\bm v$ $(\frac{\bm x_n^* + \bm x_n}{2}, \omega)\big) = 0
$ when $d = 3$. It appears that the condition $\det\big(\mathbf{D}\bm v(\frac{\bm x_n^* + \bm x_n}{2}, \omega)\big) = 0$ is highly restrictive and may not be satisfied in general cases. Conditions that are difficult to satisfy also appear when $d > 3$. To overcome this issue,  we employ the splitting method proposed in \cite{Kang1995}. We split the divergence-free vector field $\bm v(\bm x, \omega): \mathbb{R}^d \to \mathbb{R}^d$ into $d-1$ vector fields, i.e. $\bm v(\bm x, \omega) = \sum_{k=1}^{d-1} \bm v_{k}(\bm x, \omega)$, where each $\bm v_{k}(\bm x, \omega) = \bm v_{k}(x_k, x_{k+1}, \omega)$ is Hamiltonian in the variables $(x_k,x_{k+1})$. Then we can split the equation \eqref{PasTraSDEforVx} into $d$ sub-problems as following:
\begin{equation}
    \label{3DsubProbEqu}
    \left\{
    \begin{aligned}
        &\bm x^{1}(t_n) = \bm x(t_{n-1}) + \bm v_{1}(\frac{\bm x(t_{n-1}) + \bm x^{1}(t_n) }{2}, \omega) \Delta t, \\
        &\bm x^{2}(t_n) = \bm x^1(t_{n}) + \bm v_{2}(\frac{\bm x^1(t_{n}) + \bm x^{2}(t_n) }{2}, \omega) \Delta t, \\
        & \cdots \cdots \notag \\
        &\bm x^{d-1}(t_n) = \bm x^{d-2}(t_{n}) + \bm v_{d-1}(\frac{\bm x^{d-2}(t_{n}) + \bm x^{d-1}(t_n) }{2}, \omega) \Delta t, \\
        &\bm x(t_n)  = \bm x^{d-1}(t_n) + \sigma \sqrt{\Delta t} \bm w_{n-1}.
    \end{aligned}  
    \right.
\end{equation}
For example, following the Randomization Method discussed in \cite{Majda1999}, a 3D random velocity field can be generated by \eqref{vGene} with $\theta_n = 0$.
Then we can separate this velocity field into
\begin{align}
    \bm v_{1}(\bm x) =  \frac{1}{\sqrt{N}}\sum_{n=1}^N \big(\bm u_n^{1,2}(\bm k_n) \cos(\bm k_n \cdot \bm x) + \bm w_n^{1,2}(\bm k_n) \sin(\bm k_n \cdot \bm x)\big), \notag \\
    \bm v_{2}(\bm x) =  \frac{1}{\sqrt{N}}\sum_{n=1}^N \big(\bm u_n^{2,3}(\bm k_n) \cos(\bm k_n \cdot \bm x) + \bm w_n^{2,3}(\bm k_n) \sin(\bm k_n \cdot \bm x)\big) \notag,
\end{align}
here 
\begin{align}
    \bm u_n^{1,2}(\bm k_n) = (\xi_{n,2}k_{n,3} - \xi_{n,3}k_{n,2},\ \xi_{n,3}k_{n,1} - \xi_{n,2}\frac{k_{n,1} k_{n,3}}{k_{n,2}},\ 0 ), \label{u12}\\
    \bm u_n^{2,3}(\bm k_n) = (0,\ \xi_{n,2}\frac{k_{n,1} k_{n,3}}{k_{n,2}} - \xi_{n,1} k_{n,3},\ \xi_{n,1}k_{n,2} - \xi_{n,2}k_{n,1}), \label{u23} 
\end{align}
and by replacing $\xi_{n,i}$ with $\zeta_{n,i}$, where $i = 1,2,3$ in \eqref{u12}, \eqref{u23}, we get $\bm w_n^{1,2}(\bm k_n)$ and $\bm w_n^{2,3}(\bm k_n)$.
Then we adopt the first-order Lie-Trotter splitting method to discrete these three sub-problems where the subproblems are discretized by the Euler midpoint rule, separately. By
applying the splitting method in \cite{McLachlan2002}, we can split a $d$-dimensional velocity field into a summation of $d-1$ velocity fields and construct a structure-preserving scheme. More details can be found in \cite{Kang1995}, \cite{Hairer2006}. This allows us to construct the corresponding numerical scheme for each concrete realization of the $d$-dimensional random field.
\subsection{Volume-preserving property of the structure-preserving scheme }\label{Sec::CompSPandEul}
In this subsection, we would like to show the volume-preserving property of our implicit structure-preserving scheme in the stochastic sense and make a comparison with the Euler-Maruyama scheme. 
Here we use 2D random flows for illustration since, as mentioned in the previous subsection, for high dimensional random flows we are also separating them into multiple stochastic flows which are all Hamiltonian in two spatial variables. The analysis of high dimensional cases can be proved by following the idea of Theorem \ref{Thm::ErrEstforBotnSch} and Theorem \ref{Coro::subdiffforSP} in the next section. 

By following the numerical scheme 
\eqref{SubProb1Discret} \eqref{SubProb2Discret}, we can generate a solution series $\{\bm x_n^{\omega}\}, n=1,2,3,\cdots$. We can see that the solution $\bm x_{n+1}^{\omega}$ is only determined by the solution at the previous time step, $\bm x_n^{\omega}$. So it can be viewed as a Markov process.
Here we demonstrate some important properties of the random flows are preserved even after numerical discretization.
As an analogy to the continuous-time case \eqref{ContinTimeCase}, we can define the environment process as viewed from the numerical solution $\bm x_n^{\omega}$ at different discrete time steps 
\begin{equation}
    \eta_n = \tau_{\bm x_n^{\omega}} \omega(n \Delta t),\quad \eta_0 = \omega(0). 
\end{equation}
The above environment process is defined on the space of trajectories $(\Tilde{\Omega}, \mathcal{B})$ with $\Tilde{\Omega} = C([0, \infty) \cap \Delta t \mathbb{Z}; \mathbb{R}^d )$. It is a subspace of $\Omega$ with time parameter lies only on $\Delta t \mathbb{Z}$ and we still use $\mathbf{E}_{\Tilde{\bm x}}$ to denote the corresponding expectation operator. 
Then, we can also define
\begin{equation}
    \bm S_n f(\Tilde{\bm x}) = \mathbf{ME}_{\Tilde{\bm x}} f(\eta_n) .
    \label{DiscreteTimeCase}
\end{equation}

\begin{proposition}

    Under the structure-preserving numerical scheme setting, $P_0$ is an invariant measure of $\eta_n$, i.e., $P_0$ is an invariant measure of the Markov semigroup $\{ \bm S_n\}$.
    \label{P0invariantSnmp}
\end{proposition}

\begin{proposition} \label{prop42}
    Under the structure-preserving numerical scheme setting, $\bm S_n$ has the property that 
    \begin{equation}
        ||\bm S_n f||_{L^2(\mathbb{R}^d)} \leq \exp\big(-2 c_{1}(\sigma) n \Delta t\big) ||f||_{L^2(\mathbb{R}^d)}
        \label{SnExpDecayProp}
    \end{equation}
   
\end{proposition}

The proof of Proposition \ref{P0invariantSnmp} and Proposition \ref{prop42} can be found in \cite{Junlong2021}. We can obtain from Proposition \ref{P0invariantSnmp} that, in the structure-preserving scheme, $\mathbf{E}\bm S_n f = \mathbf{E}\bm S_{n-1} f$ for all $n$. 
Then we use $\bm \Theta_{\Delta t}^{\omega, S}$ to denote the numerical integrator associated with the structure-preserving scheme during $\Delta t$ time and define the transition probability density function (PDF) of the corresponding solution process,
\begin{equation}
    p^S_{\omega}(\bm x, \bm y) =   p_0(\bm \Theta_{\Delta t}^{\omega, S}(\bm x), \bm y),
    \label{Defp1omega}
\end{equation}
where $p_0$ is the kernel of the sub-process corresponding $d \bm w_t$, i.e.
\begin{equation}
    p_0(\bm x, \bm y) =  \frac{1}{2\pi \sigma^2 \Delta t}\exp \Big( - \frac{||\bm y -\bm x||^2}{2 \sigma^2 \Delta t}\Big).
    \label{Defofp0}
\end{equation}

\begin{definition}
     A numerical scheme is volume-preserving if the integral of its transition PDF is always equal to 1, i.e. 
\begin{equation}
    \int p^S_{\omega}(\bm x, \bm y) d\bm x \equiv 1.
\end{equation}
\end{definition}

Then we can find that
\begin{equation}
    \int p^S_{\omega}(\bm x, \bm y) d\bm x = \int p_0(\bm \Theta_{\Delta t}^{\omega, S}(\bm x), \bm y) d\bm x = \int p_0(\bm z,\bm y) \det(\mathbf{D}\bm \Theta_{\Delta t}^{\omega, S})^{-1} d\bm z.
    \label{InvariantMeasureP0}
\end{equation}
 Here $\bm z = \bm \Theta_{\Delta t}^{\omega, S}(\bm x)$.  It can be verified that in the structure-preserving scheme, 
\begin{equation}
    \det(\mathbf{D}\bm \Theta_{\Delta t}^{\omega,S}) = \det(\frac{\partial \bm x^*}{\partial \bm x}) \equiv 1,
\end{equation} 
always holds from \eqref{CondiforVolPre}. Together with the fact that $\int p_0(\bm z, \bm y) d\bm z = 1$, we can prove that
\begin{equation}
    \int p^S_{\omega}(\bm x, \bm y) d\bm x \equiv 1.
\end{equation}

This implies that $\bm S_n$ is volume-preserving for the structure-preserving scheme. While for the Euler-Maruyama scheme, we use $\bm \Theta_{\Delta t}^{\omega, E}$ to denote the corresponding numerical integrator and $p^E_{\omega}(\bm x, \bm y) = p_0(\bm \Theta_{\Delta t}^{\omega, E}, \bm y )$ to denote the corresponding transition PDF. To be specific
\begin{equation}
    \bm \Theta_{\Delta t}^{\omega,E} (\bm x) = \bm x + \Delta t \bm v(\bm x), 
\end{equation}
and it can be computed that 
\begin{align}
    \det(\mathbf{D}\bm \Theta_{\Delta t}^{\omega,E}) = 
    \begin{vmatrix}
        1 - \Delta t \Psi_{x_1 x_2} & - \Delta t \Psi_{x_2 x_2} \\
        \Delta t \Psi_{x_1 x_1} & 1 + \Delta t \Psi_{x_1 x_2}
    \end{vmatrix}
    =1 + (\Delta t)^2 \det(\mathbf{H}(\bm x)).
\end{align} 
Here $\mathbf{H}(\bm x)$ is the Hessian matrix of the stream function $\Psi$.


Following the randomization method in \cite{Kraichnan1970}, we can generate 2D random velocity fields in the form \eqref{vGene} with $\theta_n = 0$.
With Assumption \ref{Assforvpro} and \ref{Assforvregular}, here we consider $\bm v \in (C_b^1(\mathcal{X}))^2$, then we can find the upper bound for $|\frac{\partial v_i}{\partial x_j}|,\ i,j \in \{1, 2\} ,$ and denote it as $M_1$. Then we can obtain that

\begin{equation}
    |det(\mathbf{H}(\bm x))| = |\Psi_{x_1 x_1} \Psi_{x_2 x_2} - \Psi_{x_1 x_2}^2| \leq 2 M_1^2.
    \label{BoundofdetH}
\end{equation}
Here $\Psi$ is the stream function (Hamiltion) of the velocity field $\bm v$ and we use $\Psi_{x_i x_j}$ to denote $\frac{\partial^2 \Psi}{\partial x_i \partial x_j}$. It can be represented by the following formula, 
\begin{equation}
    \Psi (\bm x)  = \frac{1}{\sqrt{N}} \sum_{n=1}^N \xi_n \sin(\bm k_n \cdot \bm x) - \zeta_n \cos(\bm k_n \cdot \bm x) = \frac{1}{\sqrt{N}} \sum_{n=1}^N \Psi_n(\bm x),
    \label{Hamiltion}
\end{equation}
and here we define $\Psi_n (\bm x) = \xi_n \sin(\bm k_n \cdot \bm x) - \zeta_n \cos(\bm k_n \cdot \bm x)$, $n = 1, \cdots, N$. 

 It can be obtained from \eqref{BoundofdetH} that $det(\mathbf{H}(\bm x))$ is also bounded. Then we require $\Delta t < \frac{1}{2M_1}$ and then can obtain that $\frac{1}{2} \leq 1 + (\Delta t)^2 \det(\mathbf{H}(\bm x)) \leq \frac{3}{2}$. This guarantees the fraction in \eqref{P0forEulscheme} always meaningful. Then we can further compute that 
 \begin{align}
    det(\mathbf{H}(\bm x)) &=  \Psi_{x_1 x_1} \Psi_{x_2 x_2} - \Psi_{x_1 x_2}^2 \notag \\
    &= \frac{1}{N} \big((\sum_{i=1}^N k_{i,1}^2 \Psi_i(\bm x)) (\sum_{j=1}^N k_{j,2}^2 \Psi_j(\bm x)) - (\sum_{m=1}^N k_{m,1} k_{m,2} \Psi_m(\bm x))^2 \big) \notag \\
    &= \frac{1}{N}  \sum_{i > j} (k_{i,1} k_{j,2} - k_{i,2} k_{j,1})^2 \Psi_{i}(\bm x) \Psi_j(\bm x). 
 \end{align}
\begin{align}
    \Psi_i(\bm x) \Psi_j(\bm x) = & \sqrt{(\xi_i^2 +\zeta_i^2)(\xi_j^2 + \zeta_j^2)} \sin(\bm k_i \cdot \bm x + \alpha_i) \sin(\bm k_j \cdot \bm x + \alpha_j) \notag \\
    = &\frac{\sqrt{(\xi_i^2 +\zeta_i^2)\big(\xi_j^2 + \zeta_j^2)}}{2} (\cos((\bm k_i - \bm k_j) \cdot \bm x + \alpha_i - \alpha_j) \notag \\
    &- \cos((\bm k_i + \bm k_j) \cdot \bm x + \alpha_i + \alpha_j) \big).
\end{align}
Here $\alpha_1,\cdots,\alpha_N$ are obtained by the auxiliary angle formula. Due to $k_{n,1}, k_{n,2}$, $\xi_n, \zeta_n, n = 1,\cdots N$ are independent random variables, w.l.o.g we can assume that $k_{i,1} k_{j,2} - k_{i,2} k_{j,1} \neq 0$  and $\bm k_i \pm \bm k_j$ are different from each other for any $i, j$. Then the whole $det(\mathbf{H}(\bm x))$ can be treated as the summation of trigonometric functions and no terms can be eliminated from each other. So $det(\mathbf{H}(\bm x))$ is not a constant function and it oscillates up and down around 0 as $\bm x$ changes.
\begin{lemma}
\label{Lemma2forEul}
    For 2D random fields generated by \eqref{vGene}, with $\theta_n = 0, $ and $\Delta t < \frac{1}{2M_1}$ small enough, the Euler-Maruyama scheme is not a volume-preserving scheme, i.e. 
    \begin{equation}
         \int p^E_{\omega}(\bm x, \bm y) d\bm x \not \equiv 1. 
    \end{equation}
\end{lemma}
\begin{proof}


 For the Euler-Maruyama scheme, the equation \eqref{InvariantMeasureP0} will become
\begin{equation}
    \int p^E_{\omega}(\bm x, \bm y) d\bm x =  \int p_0(\bm z,\bm y) \frac{1}{1+(\Delta t)^2 det(\mathbf{H}(\bm x ))} d\bm z = \int p_0(\bm z,\bm y) h(\bm z) d\bm z.
    \label{P0forEulscheme}
\end{equation}
Here $\bm z = \bm \Theta_{\Delta t}^{\omega, E}(\bm x)$ and we use $h(\bm z)$ to denote $ \frac{1}{1+(\Delta t)^2 det(\mathbf{H}(\bm x))}$. 
Then we can find two points $\bm z_0=\bm \Theta_{\Delta t}^{\omega, E}(\bm x_0), \bm z_1=\bm \Theta_{\Delta t}^{\omega, E}(\bm x_1)$, s.t. $det(\mathbf{H}(\bm x_0)) > 0$, $  det(\mathbf{H}(\bm x_1)) < 0$, i.e. we can find
\begin{equation}
    h(\bm z_0) < 1 < h(\bm z_1),\ \bm z_0, \bm z_1 \in \mathbb{R}^2 .
\end{equation}

Due to the form of $p_0(\bm x, \bm y)$ defined in \eqref{Defofp0} and $h(\bm z)$ is continuous, the integral $ \int p_0(\bm z,\bm y) h(\bm z) d\bm z$ will be dominated by the neighborhood of $\bm y$ when $\Delta t$ is sufficient small. In this situation, by setting $\bm y_0 = \bm x_0$ and $\bm y_1 = \bm x_1$, we can obtain that 
\begin{equation}
    \int p^E_{\omega}(\bm x, \bm y_0) d\bm x \neq \int p^E_{\omega}(\bm x, \bm y_1) d\bm x
\end{equation}

So that $\int p^E_{\omega}(\bm x, \bm y) d\bm x \not \equiv 1$ in the explicit Euler-Maruyama scheme. This shows that the explicit Euler-Maruyama scheme is not volume-preserving. 
\hfill $\square$
\end{proof}

\subsection{Modified equation systems of two numerical schemes}
\label{subsec::Modifiedequsys}
In this subsection, we would like to first construct the modified equation systems of two numerical schemes that will be used in the convergence analysis. 
In the subsequent analysis, all variables depend on the implementation of random velocity fields, i.e. are dependent on $\omega$, and for writing convenience we have all ignored the $\omega$ in the writing. 
 With operator $\mathcal{L}_0$ defined in \eqref{defofL0}, we consider the backward Kolmogorov equation defined as follows,
\begin{equation}
     u_t = \mathcal{L}_0 u,\ \text{and}\ u(\bm x, 0) = \phi(\bm x). 
     \label{BackKomoEqua}
\end{equation}
 Then integrate \eqref{BackKomoEqua} from $t = 0$ to $t = \Delta t$, we can obtain that 
\begin{equation}
    u (\bm x, \Delta t) = \phi(\bm x) + \int_0^{\Delta t} \mathcal{L}_0 u(\bm x, s) ds.
\end{equation}
Then by applying weak Taylor expansion on $u(\bm x, s)$ and denoting $\frac{\partial^k}{\partial s^k} u(\bm x, 0) = \mathcal{L}_0^k \phi(\bm x)$, we can obtain that 
\begin{equation}
    u(\bm x ,\Delta t) = \phi(\bm x) + \sum_{k=0}^N \frac{\Delta t ^{k+1}}{(k+1)!} \mathcal{L}_0^{k+1} \phi(\bm x) + O(\Delta t^{N+2}).
    \label{Taylorexpanofu}
\end{equation}
For $d$-dimensional cases ($d \geq 2$), we study the flow generated by the method mentioned in subsection \ref{subsec::Constof2D} and \ref{subsec::constofnd}. We define $\mathcal{G}_{1} = \bm v_1 \cdot \nabla$, $\mathcal{G}_{2} = \bm v_2 \cdot \nabla$, ..., $\mathcal{G}_{d-1} = \bm v_{d-2} \cdot \nabla$, $\mathcal{G}_{d} =  \frac{\sigma^2}{2}  \Delta$ and compute

\begin{equation}
    \left\{
    \begin{aligned}
        &\partial_t u^1 = \mathcal{G}_{1} u^1,\ u^1(\bm x, 0) = u(\bm x, 0), \\
        &\partial_t u^2 = \mathcal{G}_{2} u^2,\ u^2(\bm x, 0) = u^1(\bm x, \Delta t), \\
        & \cdots \\
        &\partial_t u^{d-1} = \mathcal{G}_{d-1} u^{d-1},\ u^{d-1}(\bm x, 0) = u^{d-2}(\bm x, \Delta t), \\
        &\partial_t u^{d} = \mathcal{G}_{d} u^{d},\ u^{d}(\bm x, 0) = u^{d-1}(\bm x, \Delta t),\\
        & u(\bm x, 0) = \phi(\bm x).   
    \end{aligned}
    \right.
    \label{BackKomoEquofsypnd}
\end{equation}
Then $u^d (\cdot, \Delta t)$ will be the flow at time $t = \Delta t$ generated by our scheme and be an approximation to $u(\cdot, \Delta t)$.  To analyze the difference between \eqref{BackKomoEqua} and \eqref{BackKomoEquofsypnd}, we shall resort to the Baker-Campbell-Hausdorff (BCH) formula and replace the matrices with differential operators. Let $\mathcal{I}_{\Delta t} = \exp(\Delta t \mathcal{L}^{num})$ denote the composite flow operator associated with \eqref{BackKomoEquofsypnd} and $\mathcal{L}^{num}$ be the corresponding generator, i.e.,
\begin{align}
    \mathcal{I}_{\Delta t}& u(\bm x, 0) := \exp(\Delta t \mathcal{G}_{d})\exp(\Delta t \mathcal{G}_{d-1}) \cdots \exp(\Delta t \mathcal{G}_{1} ) \phi(\bm x), \notag \\
    & = \exp(\Delta t (\mathcal{G}_{d} + \cdots + \mathcal{G}_{1})) \phi(\bm x) + \frac{\Delta t^2}{2} \sum_{i > j}^{d} [\mathcal{G}_{i}, \mathcal{G}_{j}] \phi(\bm x) + O(\Delta t^3), \notag \\
    & = \exp(\Delta t \mathcal{L}_0) \phi(\bm x) + \frac{\Delta t^2}{2} \sum_{i > j}^{d} [\mathcal{G}_{i}, \mathcal{G}_{j}] \phi(\bm x) + O(\Delta t^3).
    \label{BCHformforLnum}
\end{align}

Here $[,]$ is the Lie-Bracket. Then we apply weak Taylor expansion and obtain that
\begin{equation}
     u^{d}(\bm x ,\Delta t) = \phi(\bm x) + \Delta t \mathcal{L}_0 \phi(\bm x) + \Delta t^2 \mathcal{A}_1 \phi(\bm x) + O(\Delta t^3).
     \label{Tayexpanofunum}
\end{equation}


Here $\mathcal{A}_1 = \frac{1}{2} \sum_{i > j}^{d} [\mathcal{G}_{i}, \mathcal{G}_{j}] + \frac{1}{2} \mathcal{L}_0^2$. We can also obtain an asymptotic form of the generator of this modified equation in terms of $\Delta t$
\begin{equation}
    \mathcal{L}^{num} = \mathcal{L}_0 + \Delta t \mathcal{L}_1 + \Delta t^2 \mathcal{L}_2 + \cdots.
\end{equation}

We denote the truncated generator as
\begin{equation}
    \mathcal{L}^{\Delta t, k} = \mathcal{L}_0 + \Delta t \mathcal{L}_1 + \cdots + \Delta t^k \mathcal{L}_k.
    \label{Ldeltak}
\end{equation}

Substitute \eqref{Ldeltak} into \eqref{Taylorexpanofu}. Then by comparing the corresponding term's coefficient with \eqref{Tayexpanofunum}, we can obtain that 
\begin{equation}
    \mathcal{L}_1 = \mathcal{A}_1 - \frac{1}{2} \mathcal{L}_0^2 = \frac{1}{2} \sum_{i > j}^{d} [\mathcal{G}_{i}, \mathcal{G}_{j}].
\end{equation}


Here we give a specific form for the 2-dimensional cases, while for the high-dimensional cases, there are concrete discussions in the subsequent proofs. Hence, the modified flow of $\bm x^{\Delta t, k}$ under the structure-preserving scheme can be written as
\begin{equation}
    \left\{
    \begin{aligned}
        dx_1 = &\big[-\Psi_{x_2} + \Delta t \big(\frac{1}{2} (\Psi_{x_2 x_2} \Psi_{x_1} + \Psi_{x_1 x_2} \Psi_{x_2})+ \frac{\sigma^2}{4}(\Psi_{x_1 x_1 x_2} + \Psi_{x_2 x_2 x_2}) \big)\big] dt  \\
        &+ \sigma dW_1 + \frac{\sigma}{2} \Delta t \Psi_{x_1 x_2} dW_1 + \frac{\sigma}{2} \Delta t \Psi_{x_1 x_1} dW_2, \\
        dx_2 = &\big[\Psi_{x_1} - \Delta t \big(\frac{1}{2} (\Psi_{x_1 x_2} \Psi_{x_1} + \Psi_{x_1 x_1} \Psi_{x_2})+ \frac{\sigma^2}{4}(\Psi_{x_1 x_2 x_2} + \Psi_{x_1 x_1 x_1}) \big)\big] dt  \\
        &+ \sigma dW_2 - \frac{\sigma}{2} \Delta t \Psi_{x_1 x_2} dW_2 - \frac{\sigma}{2} \Delta t \Psi_{x_2 x_2} dW_1.
    \end{aligned}
    \right.
    \label{ModifiedEquofEul}
\end{equation}

Here $\Psi$ is the stream function and we use $\Psi_{x_1}, \Psi_{x_2}$ to denote $\frac{\partial \Psi}{\partial x_1}$ and $\frac{\partial \Psi}{\partial x_2}$ respectively. Similar notations are also used in the remaining of this paper. Then we can denote 
\begin{equation}
    \Tilde{\bm v}_S := 
    \begin{pmatrix}
        \frac{1}{2} (\Psi_{x_2 x_2} \Psi_{x_1} + \Psi_{x_1 x_2} \Psi_{x_2})+ \frac{\sigma^2}{4}(\Psi_{x_1 x_1 x_2} + \Psi_{x_2 x_2 x_2}) \\
        -\frac{1}{2} (\Psi_{x_1 x_2} \Psi_{x_1} + \Psi_{x_1 x_1} \Psi_{x_2})- \frac{\sigma^2}{4}(\Psi_{x_1 x_2 x_2} + \Psi_{x_1 x_1 x_1})
    \end{pmatrix}. \notag
\end{equation}
\begin{equation}
    \Tilde{\mathbf{D}}_S := \big((\mathbf{I}_2 +\frac{ \Delta t}{2} \begin{pmatrix}
        \Psi_{x_1 x_2} &  \Psi_{x_1 x_1} \\
        -\Psi_{x_2 x_2} & - \Psi_{x_1 x_2}
    \end{pmatrix} )(\mathbf{I}_2 + \frac{\Delta t}{2} \begin{pmatrix}
        \Psi_{x_1 x_2} &  \Psi_{x_1 x_1} \\
        -\Psi_{x_2 x_2} & - \Psi_{x_1 x_2}
    \end{pmatrix})^T - \mathbf{I}_2\big) / \Delta t. \label{Dstar}
\end{equation}

Then the density function of particles $u(\bm x, t)$ obtained from the structure-preserving scheme satisfies a modified Fokker-Planck equation given by
\begin{equation}
    u_t = -\nabla((\bm v + \Delta t \Tilde{\bm v}_S)  u) + D_0 \nabla \nabla : (\mathbf{I}_2 + \Delta t\Tilde{\mathbf{D}}_S) u.
    \label{ModifiedFPEqu}
\end{equation}

The inner product between matrices is denoted by $A:B= tr(A^T B) = \sum_{i,j} a_{(i,j)}$ $ b_{(i,j)}$. It follows that $\Delta = \nabla \nabla : \mathbf{I}_2$ and $\nabla \nabla : \Tilde{\mathbf{D}}_S$ are defined accordingly.
Similar analysis can be done for the Euler-Maruyama scheme and we can obtain and denote $\Tilde{\bm v}_E$ and $\Tilde{\mathbf{D}}_E$  correspondingly, where
\begin{equation}
    \Tilde{\bm v}_E := 
    \begin{pmatrix}
        \frac{1}{2} (\Psi_{x_2 x_2} \Psi_{x_1} - \Psi_{x_1 x_2} \Psi_{x_2})+ \frac{\sigma^2}{4}(\Psi_{x_1 x_1 x_2} + \Psi_{x_2 x_2 x_2}) \\
        -\frac{1}{2} (\Psi_{x_1 x_2} \Psi_{x_1} - \Psi_{x_1 x_1} \Psi_{x_2})- \frac{\sigma^2}{4}(\Psi_{x_1 x_2 x_2} + \Psi_{x_1 x_1 x_1})
    \end{pmatrix}. \notag
\end{equation}
\begin{equation}
   \Tilde{\mathbf{D}}_E := \big((\mathbf{I}_2 +\frac{ \Delta t}{2} \begin{pmatrix}
        \Psi_{x_1 x_2} &  - \Psi_{x_1 x_1} \\
        \Psi_{x_2 x_2} & - \Psi_{x_1 x_2}
    \end{pmatrix} )(\mathbf{I}_2 + \frac{\Delta t}{2} \begin{pmatrix}
        \Psi_{x_1 x_2} &  - \Psi_{x_1 x_1} \\
        \Psi_{x_2 x_2} & - \Psi_{x_1 x_2}
    \end{pmatrix})^T - \mathbf{I}_2\big) / \Delta t.
    \label{D}
\end{equation}

\section{Convergence analysis}
In this section, we shall provide the convergence analysis of our stochastic structure-preserving scheme in computing effective diffusivity and anomalous diffusion phenomena. It has been mentioned in \cite{Fannjiang1998} that the diffusivity is always enhanced in incompressible flows when $\sigma > 0$. So the analysis of the diffusive and super-diffusive cases is straightforward, while the sub-diffusive cases are investigated with $\sigma \to 0$. 

With the modified flow system, we can observe and analyze the difference between our structure-preserving scheme and the Euler-Maruyama scheme when computing random flows with a specific energy spectrum density (Subsection \ref{subsec::ConAnaofEffDiff}). 
Meanwhile, the error caused by the structure-preserving scheme in computing the power law of random fields with super-diffusive behavior is also estimated (Subsection \ref{subsec::ErrEstofSupDiffcase}).

\subsection{Convergence analysis of two schemes in computing effective diffusivity}
\label{subsec::ConAnaofEffDiff}
In this subsection, we would like to provide the convergence analysis of two numerical schemes in computing the effective diffusivity and make a comparison between them.
To make the following analysis more concise, here we consider the 2-dimensional random flows with stream functions defined in \eqref{Hamiltion}, where $\xi_n, \zeta_n$ are i.i.d. standard Gaussian variables and $\bm k_n = (k_{n,1}, k_{n,2})$ are random variables with mean zero and finite variance. Here we require $k^2_{n,1} +k_{n,2}^2 > 0$. For 3-dimensional cases, as mentioned in subsection \ref{subsec::constofnd}, we can separate the problem into sub-problems and the velocity field in each sub-problem is Hamiltonian in two spatial variables. This is equivalent to each sub-problem being a 2-dimensional case, and then we combine these sub-problems together.

From the perspective of the stream function (Hamiltonian) $\Psi$, we can easily see the numerical diffusion error caused by the Euler-Maruyama scheme, but this is not relevant to our subsequent convergence analysis in computing effective diffusivity, so we put the analysis of this part in Appendix \ref{AppendixforStreamfct}, which is available to interested readers.

Then we denote and consider $\Phi(\bm x(t)) = \frac{||\bm x(t)||^2}{d},\ \bm x(t) = (x_1(t), \cdots, x_d(t))^T$ and the effective diffusivity of the modified flow system can be computed by $\frac{\mathbb{E} [\Phi(\bm x(t))]}{2t}$. 

\begin{remark}
    Since the random velocity fields we study are all isotropic, and thus in calculating the effective diffusivity, the value of $\frac{\mathbb{E}[\Phi(t)]}{2 t}$ is equivalent to the value of $D^{11}_E$.
\end{remark}

Due to the spectral gap property and the operator for these two numerical schemes being uniform elliptic when $\Delta t$ is small enough, we can consider for the corresponding corrector problem of these two schemes and define 
\begin{equation}
    \text{Euler-Maruyama:}\ \hat{\bm \mu}_E(t):= \int_0^{\infty} (\bm v (\bm x_E(t+s)) + \Delta t \Tilde{\bm v}_E(\bm x_E(t+s))) ds.
\end{equation}
\begin{equation}
    \text{Structure-Preserving:}\ \hat{\bm \mu}_S(t):= \int_0^{\infty} (\bm v (\bm x_S(t+s)) + \Delta t \Tilde{\bm v}_S(\bm x_S(t+s))) ds.
\end{equation}

Here we use $\bm x_E$ and $\bm x_S$ to denote the different paths of these two numerical schemes' modified equation system. And $\hat{\bm \mu}_E$ satisfies the following continuous-type corrector problem which corresponds to the modified equation system.
\begin{equation}
    \mathcal{L}^E \hat{\bm \mu}_E = - (\bm v + \Delta t \Tilde{\bm v}_E).
\end{equation}
Here $\mathcal{L}^E$ is the generator of the semigroup corresponding to the Euler-Maruyama scheme.  Then we can deduce that 
\begin{equation}
    \frac{d \hat{\bm \mu}_E}{d t} = - (\bm v (\bm x_E(t)) + \Delta t \Tilde{\bm v}_E(\bm x_E(t))) . 
\end{equation}
Similar results can also be obtained for $\hat{\bm \mu}_S$.
With $\bm x_E$ and $\bm x_S$, we can define the corresponding functions $\Phi^E$ and $\Phi^S$ for two numerical schemes, i.e.
\begin{equation}
    \Phi^E(t) = \Phi(\bm x_E(t)),\ \text{and} \  \Phi^S(t) = \Phi(\bm x_S(t)).
\end{equation}

Then we also denote $\mathbf{G}_E$ and $\mathbf{G}_S$ as the following formulas.
\begin{equation}
    \mathbf{G}_E := \sigma (\mathbf{I}_2 + \frac{\Delta t}{2} 
    \begin{pmatrix}
        \Psi_{x_1 x_2} & - \Psi_{x_1 x_1} \\
        \Psi_{x_2 x_2} & - \Psi_{x_1 x_2}
    \end{pmatrix}), \quad
    \mathbf{G}_S := \sigma (\mathbf{I}_2 + \frac{\Delta t}{2} 
    \begin{pmatrix}
        \Psi_{x_1 x_2} &  \Psi_{x_1 x_1} \\
        -\Psi_{x_2 x_2} & - \Psi_{x_1 x_2}
    \end{pmatrix}).
\end{equation}
\begin{lemma}
\label{Lemma::dPhidt}
     Let $\bm x_E(t)$ and $\bm x_S(t)$  be the solution of the modified flow system of two schemes with $\bm v$ satisfies Assumption \ref{AssforEk}. Then together with $\hat{\bm \mu}_E$, we have the following formula for $\Phi^E(t)$ in the Euler-Maruyama scheme.
     \begin{equation}
         \frac{d\Phi^E}{dt} = - \frac{d (\bm x_E \cdot \hat{\bm \mu}_E)}{dt} + (\bm v + \Delta t \Tilde{\bm v}_E)\cdot \hat{\bm  \mu}_E + \sigma^2 + O(\sigma^2 \Delta t^2) + \bm x_E\cdot \mathbf{G}_E \frac{d \bm w_t}{ dt}, 
     \end{equation}
where the integral of $\bm x_E\cdot \mathbf{G}_E \frac{d \bm w_t}{ dt}$ is a mean zero martingale, i.e. $\mathbf{E} [\int_0^t \bm x_E\cdot \mathbf{G}_E d \bm w_t]$ $ = 0$. Similar results can also be obtained for the modified flow system of the structure-preserving scheme, i.e.
\begin{equation}
    \frac{d\Phi^S}{dt} = - \frac{d (\bm x_S \cdot \hat{\bm \mu}_S)}{dt} + (\bm v + \Delta t \Tilde{\bm v}_S)\cdot \hat{\bm  \mu}_S + \sigma^2 + O(\sigma^2 \Delta t^2) + \bm x_S\cdot \mathbf{G}_S \frac{d \bm w_t}{ dt}, 
\end{equation}
\end{lemma}

\begin{theorem}
\label{Thm::ErrEstforBotnSch} 
    Let $\bm x_t,\ t>0$ be the solutions of the passive tracer model \eqref{PasTraSDEforVx} with the random field satisfies Assumption \ref{AssforEk}. Let $\hat{\bm \mu}$ be the solution of the corrector problem corresponding to \eqref{PasTraSDEforVx}. Then we have both the error of two schemes' modified equation systems in computing effective diffusivity converge to $0$ as $\Delta t \to 0$. 
    To be specific, we have 
    \begin{equation}
        \frac{\mathbb{E}[\Phi^*(t)]}{2t} = \frac{\sigma^2}{2} + \frac{\mathbb{E}[\int_0^t  \bm v(\bm x_s)\cdot \hat{\bm \mu}(s)  ds]}{2t} + \rho(\Delta t) + O(\sigma^2 \Delta t^2) + O(\frac{1}{\sqrt{t \Delta t}}).
        \label{ErrEstforEulandSP}
    \end{equation}
    Here $\rho (\Delta t) = O(\Delta t ^{\frac{2c_1}{2c_1 + c_3}})$ is a function which satisfies $\lim_{\Delta t \to 0} \rho(\Delta t) = 0$ and is independent of $t$. Equation \eqref{ErrEstforEulandSP} holds for $\Phi^*(t)$ taken as $\Phi^E(t)$ or $\Phi^S(t)$.
    Moreover, for $d$-dimensional ($d \geq 3$) cases, equation \eqref{ErrEstforEulandSP} also holds.
    
\end{theorem}
From the Assumption \ref{Assforvpro} and spectral gap condition, $\bm v$ and $\hat{\bm \mu}$ are both bounded. Then for a given time $t$, $\mathbb{E}[\int_0^t  \bm v(\bm x_s)\cdot \hat{\bm \mu}(s)  ds]$ will be a finite value.
Then for the structure-preserving scheme, we can split $\rho(\Delta t)$ into two parts, i.e. $\rho(\Delta t) = \rho_1^S(\Delta t) + \rho_2^S(\Delta t)$  and denote
\begin{equation}
    \rho_1^S(\Delta t) := \big|\int_0^t \int_0^{\infty} \langle \bm v(\bm x^S_s)\cdot \bm v(\bm x^S_{s+\tau})\rangle d\tau ds - \int_0^t \int_0^{\infty} \langle \bm v(\bm x_s)\cdot \bm v(\bm x_{s+\tau})\rangle d\tau ds\big| \Big / 2t.
    \label{Defrho1t}
\end{equation}
\begin{align}
    \rho_2^S(\Delta t) :=& \Delta t \big|\int_0^t \int_0^{\infty} \langle \bm v(\bm x^S_s) \cdot\Tilde{\bm v}_S(\bm x^S_{s+\tau})\rangle d\tau ds\ + \int_0^t \int_0^{\infty} \langle \Tilde{\bm v}_S(\bm x^S_s)\cdot \bm v(\bm x^S_{s+\tau})\rangle \notag \\
    &  d\tau ds + \Delta t \int_0^t \int_0^{\infty} \langle \Tilde{\bm v}_S(\bm x^S_s)\cdot \Tilde{\bm v}_S(\bm x^S_{s+\tau})\rangle d\tau ds\big| \Big /2t .
    \label{Defrho2t}
\end{align}

 Similar definitions can be done for the Euler-Maruyama scheme and we can obtain $\rho_1^E(\Delta t),\ \rho_2^E(\Delta t)$. Here we use $\bm x_t, \bm x^S_t$ and $\bm x^E_t$ to denote $\bm x(t), \bm x_S(t)$ and $\bm x_E(t)$. $\rho_1^S(\Delta t)$ and $\rho_1^E(\Delta t)$ will keep the original diffusion behavior, which is common to both numerical schemes. Then we consider for the difference of two numerical schemes between $\rho_2^S(\Delta t)$ and $\rho_2^E(\Delta t)$.

\begin{theorem}
\label{Coro::subdiffforSP}
     For the structure-preserving scheme, $\rho_2^S(\Delta t)$ can be bounded by the summation of one higher-order term of $\Delta t$ and one term correlated with the original convection enhanced behavior, i.e.   
    \begin{equation}
        \rho_2^S(\Delta t) = O(\Delta t \frac{\mathbb{E}[\int_0^t  \bm v(\bm x_s)\cdot \hat{\bm \mu}(s)  ds]}{2t})  + O(\Delta t^{\varrho}).
    \end{equation}
    Here $\varrho > 1$.
    Moreover, when $\sigma^2$ has the same scale as a fixed $\Delta t$, $\rho_2^S(\Delta t)$ in the structure-preserving method can keep the original sub-diffusion behavior of the system when $\sigma$ vanishes, i.e.
    \begin{equation}
        \lim_{t\to \infty} \frac{ \rho_2^S(\Delta t) }{\frac{\sigma^2}{2} + \frac{\mathbb{E}[\int_0^t  \bm v(\bm x_s)\cdot \hat{\bm \mu}(s)  ds]}{2t}} 
        = O(\Delta t^{\varrho - 1}) < \infty
   \end{equation}
 While the Euler-Maruyama method does not hold this property.
\end{theorem}

The proof of Lemma \ref{Lemma::dPhidt}, Theorem \ref{Thm::ErrEstforBotnSch} and Theorem \ref{Coro::subdiffforSP} can be seen in the subsection \ref{subsec::proofofThm}. Theorem \ref{Coro::subdiffforSP} and the following Corollary \ref{Coro::diffSP} also holds for 3-dimensional cases, the detailed explanation we put together in the proof of Theorem \ref{Coro::subdiffforSP}. 

\begin{corollary}
\label{Coro::diffSP}
    When diffusion coefficient $\sigma^2$ has the same scale as a fixed $\Delta t$, the error terms $\rho_2^S(\Delta t)$ caused by the structure-preserving will keep the behavior 
   of the original random velocity field with energy spectral $E(k)$ which will perform the diffusion behavior, i.e.
    \begin{equation}
        \lim_{t\to \infty} \frac{\rho_2^S(\Delta t)}{\frac{\sigma^2}{2} + \frac{\mathbb{E}[\int_0^t  \bm v(\bm x_s)\cdot \hat{\bm \mu}(s)  ds]}{2t}} 
        < O(\Delta t) + O(\Delta t^{\varrho - 1}) < \infty.
   \end{equation}
   
\end{corollary}
\begin{proof}
    for the field with diffusion behavior, the term $\lim_{t\to \infty}$ $ \frac{\mathbb{E}[\int_0^t  \bm v(\bm x_s)\cdot \hat{\bm \mu}(s)  ds]}{2t}$ will be a non-zero constant and is much larger than $\sigma^2$ when $\sigma^2$ has the same scale as a fixed $\Delta t$. Then we can obtain that
    \begin{align}
        \lim_{t\to \infty} \frac{\rho_2^S(\Delta t)}{\frac{\sigma^2}{2} + \frac{\mathbb{E}[\int_0^t  \bm v(\bm x_s)\cdot \hat{\bm \mu}(s)  ds]}{2t}} 
        &= \lim_{t\to \infty} \frac{h_1\Delta t \frac{\mathbb{E}[\int_0^t  \bm v(\bm x_s)\cdot \hat{\bm \mu}(s)  ds]}{2t} + h_2\Delta t^{\varrho}}{h_3\Delta t + \frac{\mathbb{E}[\int_0^t  \bm v(\bm x_s)\cdot \hat{\bm \mu}(s)  ds]}{2t} } \notag \\
        &< O(\Delta t) + O(\Delta t^{\varrho - 1}) < \infty. 
   \end{align}
   \hfill $\square$
\end{proof}

\subsection{Error estimate for structure-preserving scheme in super-diffusive cases}
\label{subsec::ErrEstofSupDiffcase}
For random fields with energy spectral $E(k)$ which satisfies \eqref{sharpCondi2D}, we can define the corresponding corrector problem well due to the spectral gap condition \eqref{SnExpDecayProp}. Instead, for those random fields with energy spectral that does not satisfy condition \eqref{sharpCondi2D} and exhibit super diffusion performance, $\hat{\bm \mu}_{E}$ and $\hat{\bm \mu}_{S}$ may diverge and can not be well defined. 
\begin{proposition}
\label{PropofFann}
    For 2D isotropic random flows with velocity spectrum 
    \begin{equation}
        \hat{R}_{i,j}(\bm k) \sim \frac{1}{|\bm k|^{2\alpha}} (\delta_{i,j} - \frac{k_i k_j}{|\bm k|^2}),
        \label{FannVeloSpec}
    \end{equation}
    here $\ 0 < \alpha < 1,\ |\bm k| \leq L$. Then it can be obtained that the dispersion of this system will have the following scaling law after $t$ is large enough, 
    \begin{equation}
        \mathbb{E}[(\bm x_t)^2] \sim t^{\frac{2}{2 - \alpha}}.
    \end{equation}
\end{proposition}

The proof of Proposition \ref{PropofFann} can be found in \cite{Fannjiang1998} , then the two-point correlation function $R_{i,j}(\bm x)$ is asymptotically 
$R_{i,j}(\bm x) \sim \frac{1}{|\bm x|^{2(1-\alpha)}}, \ |\bm x| \gg 1$. The random velocity field $\bm v$ which has the velocity spectrum mentioned in Proposition \ref{PropofFann} and its corresponding stream function can be realized by 
\begin{equation}
    \Psi(\bm x) = \frac{1}{\sqrt{N}} \sum_{i=1}^N\frac{1}{|\bm k_i|} (\xi_i \sin(\bm k_i \cdot \bm x) - \zeta_i\cos(\bm k_i \cdot \bm x)),\ and\ \bm v(\bm x) = \nabla^{\perp} \Psi(\bm x).
\end{equation}

Here $\forall\ 1\leq i \leq N$, $\xi_i, \zeta_i$ are i.i.d standard Gaussian random variables, the direction of $\bm k_i$ is uniformly distributed on a unit circle and the length of $\bm k_i$ satisfies the density function $\rho(|\bm k_i|)\sim 1/|\bm k_i|^{2\alpha -1}, 0 \leq |\bm k_i| \leq L$. Here the form of $\Psi(\bm x)$ is similar but different from \eqref{Hamiltion}. Then we have the following Theorem \ref{ThmforSuperdiff} for the structure-preserving scheme and its proof can be seen from subsection \ref{subsec::proofofThm}.
\begin{theorem}
    \label{ThmforSuperdiff}
    Let $\bm x_t^S$ be the solution of the modified flow system of the structure-preserving scheme at time $t$ with the original random field possessing an isotropic velocity spectrum satisfying \eqref{FannVeloSpec}. Then for $\frac{1}{4} \leq \alpha < 1$ and $t$ large enough, the modified flow system will exhibit super-diffusion performance, and the scaling law of its dispersion is consistent with the original system, i.e.
    \begin{equation}
        \mathbf{E}[(\bm x_t^S)^2] \sim t^{\frac{2}{2 - \alpha}},\ \frac{1}{4} \leq \alpha < 1.
    \end{equation}
\end{theorem}

\subsection{Proof of Lemma and Theorems}
\label{subsec::proofofThm}

\paragraph{Proof of Lemma \ref{Lemma::dPhidt}}
\begin{proof}
By applying Ito's formula on $\Phi^E(\bm x(t))$ of the modified system of the Euler-Maruyama scheme, we can obtain that
\begin{equation}
    d\Phi^E = (\bm x_E\cdot (\bm v + \Delta t \Tilde{\bm v}_E) + \frac{\sigma^2}{2} \text{Tr}(\mathbf{G}_E^T \mathbf{G}_E)) dt + \bm x_E\cdot \mathbf{G}_E d\bm w_t.
\end{equation}

Then we can compute that 
\begin{equation}
    \frac{\sigma^2}{2} \text{Tr}(\mathbf{G}_E^T \mathbf{G}_E) = \sigma^2 + \frac{\sigma^2 \Delta t^2}{8} (\Psi_{x_1 x_1}^2  + 2 \Psi_{x_1 x_2}^2 + \Psi_{x_2 x_2}^2),
\end{equation}
and denote the following integral as
\begin{align}
    \rho_E &:= \int_0^t \bm x_E(s) \cdot (\bm v (\bm x_E(s)) + \Delta t \Tilde{\bm v}_E(\bm x_E(s))) ds \notag \\
    &=  \int_0^t (\bm v + \Delta t \Tilde{\bm v}_E)\cdot \hat{\bm \mu}_E ds- [\bm x_E \cdot \hat{\bm \mu}_E]_0^t. 
\end{align}
A similar definition can also be obtained for $\rho_S$. 
\begin{align}
    \rho_S &:= \int_0^t \bm x_S(s) \cdot (\bm v (\bm x_S(s)) + \Delta t \Tilde{\bm v}_S(\bm x_S(s))) ds \notag \\
    &=  \int_0^t (\bm v + \Delta t \Tilde{\bm v}_S)\cdot \hat{\bm \mu}_S ds- [\bm x_S \cdot \hat{\bm \mu}_S]_0^t.
\end{align}

So consider the derivative to time $t$ and move the Brownian motion related term inside the martingale related term, we can get
\begin{equation}
   \frac{d\Phi^E}{dt} = - \frac{d(\bm x_E \cdot \hat{\bm \mu}_E)}{dt} + (\bm v + \Delta t \Tilde{\bm v}_E)\cdot \hat{\bm  \mu}_E + \sigma^2 + O(\sigma^2 \Delta t^2) + \bm x_E \cdot \mathbf{G}_E \frac{d \bm w_t}{ dt}
   \label{dPhiofEul}
\end{equation}
where $\mathbb{E} [\int_0^t \bm x_E \cdot \mathbf{G}_E d \bm w_t] = 0$. For the structure-preserving scheme, we can similarly obtain that 
\begin{equation}
    \frac{d\Phi^S}{dt} = - \frac{d(\bm x_S \cdot \hat{\bm \mu}_S)}{dt} + (\bm v + \Delta t \Tilde{\bm v}_S)\cdot \hat{\bm  \mu}_S + \sigma^2 + O(\sigma^2 \Delta t^2) + \bm x_S \cdot \mathbf{G}_S \frac{d \bm w_t}{ dt}, 
\end{equation}
here $\mathbb{E} [\int_0^t \bm x_S \cdot \mathbf{G}_S d \bm w_t] = 0$.
\hfill $\square$
\end{proof}

\paragraph{Proof of Theorem \ref{Thm::ErrEstforBotnSch}}

\begin{proof}
For the Euler-Maruyama scheme, we can obtain from Lemma \ref{Lemma::dPhidt} that 
\begin{equation}
         \frac{d\Phi^E}{dt} = - \frac{d(\bm x_E \cdot \hat{\bm \mu}_E)}{dt} + (\bm v + \Delta t \Tilde{\bm v}_E)\cdot \hat{\bm  \mu}_E + \sigma^2 + O(\sigma^2 \Delta t^2) + \bm x_E\cdot \mathbf{G}_E \frac{d \bm w_t}{ dt}, 
\end{equation}
It has been proved in \cite{Junlong2021} that the term $- [\bm x_E \cdot \hat{\bm \mu}_E]_0^t / 2t$ is $O(\frac{1}{t \sqrt{\Delta t}})$. It is negligible  when $t$ is large enough. And then we consider the term $(\bm v + \Delta t \Tilde{\bm v}_E)\cdot \hat{\bm  \mu}_E$ at time $t$. We want to first prove that the solution $\hat{\bm \mu}_E$ converges to the solution of the original system's corrector problem $\hat{\bm \mu}$ as $\Delta t \to 0$.

Considering the spectral gap condition, we can find a truncation time $T_0 = N\Delta t$ and obtain the following inequalities
\begin{equation}
    ||\int_{T_0}^{\infty} (\bm v(\bm x_E(s)) + \Delta t \Tilde{\bm v}_E(\bm x_E(s)) ) ds||_{L^2} \leq \frac{1}{2c_1} \exp(-2c_1 T_0),
\end{equation}
\begin{equation}
    ||\int_{T_0}^{\infty} \bm v(\bm x(s))   ds||_{L^2} \leq \frac{1}{2c_1} \exp(-2c_1 T_0).
\end{equation}

Then for any $\epsilon > 0$, we can set the truncation time $T_0$ big enough to let the following inequality hold
\begin{equation}
    ||\int_{T_0}^{\infty} (\bm v(\bm x_E(s)) + \Delta t \Tilde{\bm v}_E(\bm x_E(s)) ) ds - \int_{T_0}^{\infty} \bm v(\bm x(s))   ds||_{L^2} \leq \frac{1}{c_1} \exp(-c_1 T_0) < \epsilon.
\end{equation}

Then we consider the error before the truncation time $T_0$
\begin{align}
    ||\int_{0}^{T_0}& (\bm v(\bm x_E(s)) + \Delta t \Tilde{\bm v}_E(\bm x_E(s)) ) ds - \int_{0}^{T_0} \bm v(\bm x(s))   ds||_{L^2} \notag \\
    &\leq ||\int_{0}^{T_0} \bm v(\bm x_E(s)) ds - \int_{0}^{T_0} \bm v(\bm x(s))   ds||_{L^2} + || \int_0^{T_0}\Delta t \Tilde{\bm v}_E(\bm x_E(s)) ds||_{L^2} \notag \\
    &= ||\int_0^{T_0} (\bm v(\bm x_E(s)) - \bm v(\bm x(s))) ds||_{L^2} + \Delta t || \int_0^{T_0} \Tilde{\bm v}_E(\bm x_E(s)) ds||_{L^2}.
\end{align}

From Assumption \ref{Assforvpro} that $\bm v$ is bounded, we can obtain that $\Delta t \Tilde{\bm v}_E$ is also bounded. Then we can obtain that $\Delta t || \int_0^{T_0} \Tilde{\bm v}_E(\bm x_E(s)) ds||_{L^2}$ is $O(T_0 \Delta t) = O(N \Delta t^2)$. Also from the error estimate for Euler-Maruyama scheme in \cite{Kloeden1992}, we know that $||\int_0^{T_0} (\bm v(\bm x_E(s)) - \bm v(\bm x(s))) ds||_{L^2}$ is $O(T_0 \Delta t) = O(N\Delta t^2)$. Then we can combine these two terms together and obtain that 
\begin{equation}
    ||\int_{0}^{T_0} (\bm v(\bm x_E(s)) + \Delta t \Tilde{\bm v}_E(\bm x_E(s)) ) ds - \int_{0}^{T_0} \bm v(\bm x(s))   ds||_{L^2} \leq c_2 T_0 \Delta t.
    \label{Estof0toT0}
\end{equation}
here $c_2$ is a constant. Then we can choose $\Delta t < \epsilon / (c_2 T_0)$. And we can choose $\epsilon \approx \frac{1}{c_1} \Delta t^{\frac{2c_1}{2c_1 + c_3}}$ to balance each value of $\epsilon$ where $c_2 = \exp(c_3 T_0)$. 

Now we finish the proof that the solution $\hat{\bm \mu}_E$ converges to the solution of the original system's corrector problem $\hat{\bm \mu}$ as $\Delta t \to 0$. A similar analysis can be done for the solution $\hat{\bm \mu}_S$. Then we go back to the equation 
\eqref{dPhiofEul}. We can obtain that the long-time behavior of $(\bm v + \Delta t \Tilde{\bm v}_E) \cdot \hat{\bm  \mu}_E$ will converge to the original system's behavior as $\Delta t \to 0$. 
Also, we can further get
\begin{align}
    \frac{\mathbb{E}[\int_0^t (\bm v(\bm x_E(s)) + \Delta t \Tilde{\bm v}_E (\bm x_E(s))) \cdot \hat{\bm \mu}_E(s) ds ]}{2t} = \frac{\mathbb{E}[\int_0^t \bm v(\bm x(s))\cdot \hat{\bm \mu}(s) ds ]}{2t} + \rho(\Delta t).
\end{align}  
Here $\rho (\Delta t) = O(\Delta t ^{\frac{2c_1}{2c_1 + c_3}})$ is a function which satisfies $\lim_{\Delta t \to 0} \rho(\Delta t) = 0$ and is independent of $t$.

Then we consider d-dimensional cases ($d \geq 3$) processed under the method mentioned in subsection \ref{subsec::constofnd}. 
If we separate the velocity field $\bm v(x_1, x_2, \cdots, x_d) = \sum_{i=1}^{d-1}\bm v_i(x_i, x_{i+1})$ and use the Euler-Maruyama scheme in the advection sub-problems, we can compute and obtain the following modified equation system

\begin{equation}
    d\bm x_E = \big(\sum_{i=1}^{d-1} \bm v_i + \Delta t \sum_{j=1}^{d-1}\Tilde{\bm v}_j^E\big)  dt + (\sigma \mathbf{I}_d + \Delta t \Tilde{\mathbf{D}}_E) d\bm w_t 
\end{equation}
Here 
\begin{equation}
    \Tilde{\bm v}_j^E = 
    \begin{pmatrix} 
        0  \\
        \vdots \\
        0 \\
    
        \frac{1}{2} (\Psi^{(j)}_{x_{j+1} x_{j+1}} \Psi^{(j)}_{x_j} - \Psi^{(j)}_{x_j x_{j+1}} \Psi^{(j)}_{x_{j+1}})+ \frac{\sigma^2}{4}(\Psi^{(j)}_{x_j x_j x_{j+1}} + \Psi^{(j)}_{x_{j+1} x_{j+1} x_{j+1}}) \\
        -\frac{1}{2} (\Psi^{(j)}_{x_j x_{j+1}} \Psi^{(j)}_{x_j} - \Psi^{(j)}_{x_j x_j} \Psi^{(j)}_{x_{j+1}})- \frac{\sigma^2}{4}(\Psi^{(j)}_{x_j x_{j+1} x_{j+1}} + \Psi^{(j)}_{x_j x_j x_j}) \\
        0\\
        \vdots \\
        0 \\
    \end{pmatrix}. \notag
\end{equation}


Here $\Psi^{(j)}$ is the Hamiltonian of $\bm v_j$ in its corresponding two spatial variables $x_j, x_{j+1}$. The first $j-1$ components and the last $d-j-1$ components of $\Tilde{\bm v}_j^E$  are all $0$. From the previous analysis, the additional coefficient of the Brownian motion term, $\Delta t \Tilde{\mathbf{D}}_E$, will result in an error term that can be included by $O(\sigma^2 \Delta t)$ in computing the effective diffusivity. Then we just treat $\sum_{i=1}^{d-1} \bm v_i = \bm v$ and $\Delta t \sum_{j=1}^{d-1} \Tilde{\bm v}_j^E  = \Delta t \Tilde{\bm v}_E$ and then applying Ito's formula on $\Phi^E(t)$. The subsequent proof is independent of dimensionality and, following the idea of the previous proof, we obtain that \eqref{ErrEstforEulandSP} also holds for the Euler-Maruyama scheme in $d$-dimensional cases. 

If we use the structure-preserving scheme in the advection sub-problems, we can compute and obtain the following modified equation system 

\begin{equation}
   d\bm x_S = \big(\sum_{i=1}^{d-1} \bm v_i + \Delta t \sum_{j=1}^{d-1}\Tilde{\bm v}_j^S\big)  dt + (\sigma \mathbf{I}_d + \Delta t \Tilde{\mathbf{D}}_S) d\bm w_t 
\end{equation}
Here 
\begin{equation}
    \Tilde{\bm v}_j^S = 
    \begin{pmatrix} 
        0  \\
        \vdots \\
        0 \\
    
        \frac{1}{2} (\Psi^{(j)}_{x_{j+1} x_{j+1}} \Psi^{(j)}_{x_j} + \Psi^{(j)}_{x_j x_{j+1}} \Psi^{(j)}_{x_{j+1}})+ \frac{\sigma^2}{4}(\Psi^{(j)}_{x_j x_j x_{j+1}} + \Psi^{(j)}_{x_{j+1} x_{j+1} x_{j+1}}) \\
        -\frac{1}{2} (\Psi^{(j)}_{x_j x_{j+1}} \Psi^{(j)}_{x_j} + \Psi^{(j)}_{x_j x_j} \Psi^{(j)}_{x_{j+1}})- \frac{\sigma^2}{4}(\Psi^{(j)}_{x_j x_{j+1} x_{j+1}} + \Psi^{(j)}_{x_j x_j x_j}) \\
        0\\
        \vdots \\
        0 \\
    \end{pmatrix}. \notag
\end{equation}

And the first $j-1$ components and the last $d-j-1$ components of $\Tilde{\bm v}_j^S$ are all $0$. Then a similar analysis can be done for the structure-preserving scheme.

\hfill $\square$
\end{proof}

\paragraph{Proof of Theorem \ref{Coro::subdiffforSP}}

\begin{proof}
Consider for the structure-preserving scheme when diffusion coefficient $\sigma^2$ has the same scale as a fixed $\Delta t$, we can obtain that
\begin{align}
    (\bm v + &\Delta t \Tilde{\bm v}_S)\cdot \hat{\bm \mu}_S= \bm v\cdot \int_0^{\infty} \bm v(\bm y(s)) ds + \Delta t (\Tilde{\bm v}_S \cdot \int_0^{\infty}\bm  v(\bm y(s)) ds +  \notag \\
    & \bm v\cdot \int_0^{\infty} \Tilde{\bm v}_S(\bm y(s)) ds)+ \Delta t^2 \Tilde{\bm v}_S \cdot \int_0^{\infty} \Tilde{\bm v}_S(\bm y(s)) ds, \ \text{where}\ \bm y(s) = \bm x_S(t + s).
\end{align}

Here we use $\bm x^S_t$ to denote $\bm x_S(t)$. Then we can obtain that, 
\begin{align}
    \mathbb{E}&[\Phi(t)] =  \int_0^t \int_0^{\infty} \langle \bm v(\bm x^S_s)\cdot \bm v(\bm x^S_{s+\tau})\rangle d\tau ds  + \Delta t \int_0^t \int_0^{\infty} \langle \bm v(\bm x^S_s) \cdot\Tilde{\bm v}_S(\bm x^S_{s+\tau})\rangle d\tau ds\ \notag \\
    & +\Delta t\int_0^t \int_0^{\infty} \langle \Tilde{\bm v}_S(\bm x^S_s)\cdot \bm v(\bm x^S_{s+\tau})\rangle d\tau ds  + \Delta t^2 \int_0^t \int_0^{\infty} \langle \Tilde{\bm v}_S(\bm x^S_s)\cdot \Tilde{\bm v}_S(\bm x^S_{s+\tau})\rangle d\tau ds\notag \\
    &-\langle [\bm x^{S}\cdot \hat{\bm \mu}_S]_0^t \rangle + \sigma^2 t + O(\sigma^2 \Delta t^2 t).
    \label{EPhiforModifiedAna}
\end{align}

Then for the energy spectrum $E(k)$ with $\sigma \to 0$, the corresponding terms for the structure-preserving scheme can be represented by 
\begin{align}
    |\int_0^t &\int_0^{\infty} \langle \bm v(\bm x^S_s) \cdot \Tilde{\bm v}_S(\bm x^S_{s+\tau})\rangle d\tau ds|
    \notag \\
    <&\ \frac{1}{2}| \int_0^t \int_0^{\infty} \langle \Psi_{x_1 x_2}(\bm v(\bm x^S_s)\cdot \bm v(\bm x^S_{s+\tau}))\rangle d\tau ds| \notag \\
    & +\frac{r^2 }{2} |\int_0^t \int_0^{\infty} \langle \Psi\ (v_1(\bm x^S_s) v_2(\bm x^S_{s+\tau}))\rangle d\tau ds| \notag \\
    & + \frac{\sigma^2 r^2}{4}  |\int_0^t \int_0^{\infty} \langle \bm v(\bm x^S_s) \cdot \Tilde{\bm v}_S(\bm x^S_{s+\tau})\rangle d\tau ds|\notag \\
    \leq &\ (\frac{C_1 }{2} + \frac{\sigma^2 r^2 }{4}) |\int_0^t \int_0^{\infty} \langle \bm v(\bm x^S_s)\cdot \bm v(\bm x^S_{s+\tau})\rangle d\tau ds| \notag \\
    = &\ (\frac{C_1 }{2} + \frac{\sigma^2 r^2 }{4}) \mathbb{E}[\int_0^t \bm v(\bm x^S_s) \cdot \hat{\bm \mu}(\bm x^S_s)   ds].
    \label{SPErrEst}
\end{align}

Here $\hat{\bm \mu}(\bm x^S_s) = \int_0^{\infty} \bm v(\bm x^S_{\tau + s}) d \tau$. Then we consider the following error estimate 
\begin{equation}
    ||\int_0^t \bm v(\bm x^S_s) \cdot \hat{\bm \mu}(\bm x^S_s) ds - \int_0^t \bm v(\bm x_s) \cdot \hat{\bm \mu}(s) ds||_{L^2} \leq C_2 t (\Delta t)^{\frac{2c_1}{2c_1+c_3}}.
    \label{ErrEstofvmus}
\end{equation}

Here $C_1$ and $C_2$ are two positive constants. Now we substitute \eqref{ErrEstofvmus} back \eqref{SPErrEst} and obtain that
\begin{align}
    |\Delta t\int_0^t \int_0^{\infty} \langle \bm v(\bm x^S_s) \cdot \Tilde{\bm v}_S(\bm x^S_{s+\tau})\rangle d\tau ds| \leq &(\frac{C_1}{2} + \frac{r^2\sigma^2 }{4})  \Delta t \mathbb{E}[\int_0^t \bm v(\bm x_s) \cdot \hat{\bm \mu}(s)  ds] \notag \\
    & + (\frac{C_1 C_2 }{2} + \frac{C_2 r^2\sigma^2 }{4}) t \Delta t^{\varrho} .
    \label{Errestforvvstartotal}
\end{align}
here $\varrho = 1+ \frac{2c_1}{2c_1+c_3} >1$. A similar analysis can be done for the third term on the right-hand side of \eqref{EPhiforModifiedAna}. The fourth term on the right-hand side of \eqref{EPhiforModifiedAna} is a high-order term of $\Delta t$.  For d-dimensional cases using the structure-preserving scheme, the term $\int_0^t \int_0^{\infty} \langle \bm v(\bm x^S_s) \cdot \Tilde{\bm v}_S(\bm x^S_{s+\tau})\rangle d\tau ds$ can be transformed into

\begin{align}
    &\int_0^t \int_0^{\infty} \langle \bm v(\bm x^S_s) \cdot \Tilde{\bm v}_S(\bm x^S_{s+\tau})\rangle d\tau ds = \int_0^t \int_0^{\infty} \langle (\sum_{i=1}^{d-1} \bm v_i(\bm x_s^S)) \cdot (\sum_{j=1}^{d-1}\Tilde{\bm v}^S_j (\bm x^S_{s+\tau}))\rangle d\tau ds \notag \\
    &= \int_0^t \int_0^{\infty} \Big(\sum_{i=1}^{d-1}\langle \bm v_i(\bm x^S_s) \cdot \Tilde{\bm v}^S_i(\bm x^S_{s+\tau})\rangle  + \sum_{|i - j| = 1}\langle \bm v_i(\bm x^S_s) \cdot \Tilde{\bm v}^S_j(\bm x^S_{s+\tau})\rangle \Big) d\tau ds.
    \label{vvstarin3D}
\end{align}

The analysis of the terms $\int_0^t \int_0^{\infty} \langle \bm v_i(\bm x^S_s) \cdot \Tilde{\bm v}^S_i(\bm x^S_{s+\tau})\rangle d\tau ds$ can just follow the idea of \eqref{SPErrEst} and \eqref{ErrEstofvmus}. Then for the other terms in \eqref{vvstarin3D}, we have
\begin{align}
    &\int_0^t \int_0^{\infty} \langle \bm v_i(\bm x^S_s)  \cdot \Tilde{\bm v}^S_j (\bm x^S_{s+\tau}) \rangle d\tau ds \notag \\ 
    &= \int_0^t \int_0^{\infty} \Big\langle \frac{1}{2} \Psi^{(i)}_{x_i}(\bm x^S_s)\Big(\Psi^{(j)}_{x_{j+1} x_{j+1}}(\bm x^S_{s+\tau}) \Psi^{(j)}_{x_j}(\bm x^S_{s+\tau}) + \Psi^{(j)}_{x_j x_{j+1}}(\bm x^S_{s+\tau}) \cdot \notag \\
    &\ \Psi^{(j)}_{x_{j+1}}(\bm x^S_{s+\tau})\Big)+ \frac{\sigma^2}{4} \Psi^{(i)}_{x_i}(\bm x^S_{s}) \Big(\Psi^{(j)}_{x_j x_j x_{j+1}} (\bm x^S_{s+\tau})+ \Psi^{(j)}_{x_{j+1} x_{j+1} x_{j+1}}(\bm x^S_{s+\tau}) \Big)  \Big\rangle d\tau ds \notag \\
    & = 0.
\end{align}

This is due to the random variables in $\Psi^{(i)}$ and $\Psi^{(j)}$ are independent and mean zero. So equation \eqref{Errestforvvstartotal} also holds in $d$-dimensional cases. 
Then we can obtain that for the structure-preserving scheme
\begin{align}
    \frac{\mathbb{E}[\Phi(t)]}{2t} = &\frac{\sigma^2}{2} + \frac{\mathbb{E}[\int_0^t  \bm v(\bm x_s)\cdot \hat{\bm \mu}(s)  ds]}{2t} + \rho_1(\Delta t) + O(\Delta t \frac{\mathbb{E}[\int_0^t  \bm v(\bm x_s)\cdot \hat{\bm \mu}(s)  ds]}{2t}) \notag \\
    & + O(\Delta t^{\varrho}) + O(\sigma^2 \Delta t^2) + O(\frac{1}{\sqrt{t \Delta t}}).
\end{align}

For the structure-preserving scheme, the original error estimate term $\rho(\Delta t)$ can be separated into $\rho_1^S(\Delta t) + O(\Delta t \frac{\mathbb{E}[\int_0^t  \bm v(\bm x_s)\cdot \hat{\bm \mu}(s)  ds]}{2t}) + O(\Delta t^{\varrho})$. Here $\rho_1^S(\Delta t)$ is just the error estimate for \eqref{ErrEstofvmus}, $O(\Delta t ^{\frac{2c_1}{2c_1 + c_3}})$, which is caused by the original random velocity field under different paths. The other terms caused by the modified velocity field can be represented by $O(\Delta t \frac{\mathbb{E}[\int_0^t  \bm v(\bm x_s)\cdot \hat{\bm \mu}(s)  ds]}{2t}) + O(\Delta t^{\varrho})$. This suggests that the two additional errors introduced in the structure-preserving scheme are either higher-order terms of $\Delta t$ or correlated with the original convection-enhanced behavior. Then when $\sigma^2$ has the same scale as a fixed $\Delta t$, the other two error terms of the structure-preserving method for the original system with sub-diffusion behavior when $\sigma^2$ has the same scale as a fixed $\Delta t$ will have the following analysis
    \begin{align}
        \lim_{t\to \infty} \frac{ \rho_2^S(\Delta t)}{\frac{\sigma^2}{2} + \frac{\mathbb{E}[\int_0^t  \bm v(\bm x_s)\cdot \hat{\bm \mu}(s)  ds]}{2t}} 
        &= \lim_{t\to \infty} \frac{h_1 \Delta t \frac{\mathbb{E}[\int_0^t  \bm v(\bm x_s)\cdot \hat{\bm \mu}(s)  ds]}{2t} + h_2 \Delta t^{\varrho}}{h_3 \Delta t } \notag \\
        &= 0 + O(\Delta t^{\varrho - 1}) < \infty.
   \end{align}

Here $h_1, h_2$, and $h_3$ are all finite constants. While for the Euler-Maruyama scheme, it does not hold this behavior. For $\rho_2^E(\Delta t)$ in a 2-dimensional case, it will become that
\begin{align}
    \int_0^t &\int_0^{\infty} \langle \bm v(\bm x^E_s) \cdot \Tilde{\bm v}_E(\bm x^E_{s+\tau})\rangle d\tau ds = \frac{1}{2} \int_0^t \int_0^{\infty} \langle \Psi_{x_1 x_2}(\bm v(\bm x^E_s) \cdot \mathbf{D}_{2}\bm v(\bm x^E_{s+\tau}))\rangle d\tau ds
    \notag \\
    & + \frac{1 }{2} \int_0^t \int_0^{\infty} \langle (\Psi_{x_1 x_1} - \Psi_{x_2 x_2}) v_1(\bm x^E_s) v_2(\bm x^E_{s+\tau}))\rangle d\tau ds .
\end{align}
Here 
\begin{equation}
    \mathbf{D}_{2} = \begin{pmatrix}
        1 & 0 \\
        0 & -1
    \end{pmatrix}.
\end{equation}

It can be seen here that the errors introduced by the modified system corresponding to the Euler-Maruyama scheme cannot be transformed into a form that is consistent with the original random velocity field's convection-enhanced behavior.
\hfill $\square$
\end{proof}

\paragraph{Proof of Theorem \ref{ThmforSuperdiff}}
\begin{proof}
    
 We first compute the stream function for the structure-preserving scheme's modified equation system
\begin{equation}
    \Psi^s = \Psi - \frac{\Delta t}{2} \Psi_{x_1} \Psi_{x_2} - \frac{\sigma^2}{4} \Delta \Psi.
\end{equation}
Then we consider the correlation function of the modified system,
\begin{equation}
    R_{ij}^s(\bm r) = \langle \Psi_{x_1}(\bm r) \Psi_{x_1}(0) + \Psi_{x_2}(\bm r) \Psi_{x_2}(0)\rangle + \frac{\sigma^4 \Delta t^2}{16}\Tilde{R_1}(\bm r) + \frac{\Delta t^2}{4} \Tilde{R_2}(\bm r) - \frac{\sigma^2\Delta t}{4} \Tilde{R_3}(\bm r).
\end{equation}
Here 
\begin{align}
    \Tilde{R_1}(\bm r) 
     = \Big\langle & \sum_{i \neq j} \Big(\big(\Psi_{x_i x_i x_j}(\bm r) + \Psi_{x_j x_j x_j}(\bm r) \big) \big(\Psi_{x_i x_i x_j}(0) + \Psi_{x_j x_j x_j}(0) \big)\Big)\Big\rangle. \notag \\ \\
    \Tilde{R_2}(\bm r)  
     = \Big\langle & \sum_{i \neq j} \Big(\big(\Psi_{x_i}(\bm r) \Psi_{x_j x_j}(\bm r) + \Psi_{x_j}(\bm r) \Psi_{x_i x_j}(\bm r) \big) \cdot \notag \\
     & \quad \quad \quad \big(\Psi_{x_i}(0) \Psi_{x_j x_j}(0) + \Psi_{x_j}(0) \Psi_{x_i x_j}(0) \big)\Big)\Big\rangle. \notag \\ \\
    \Tilde{R_3}(\bm r) 
    =\Big\langle & \sum_{i \neq j} \Big(\Psi_{x_i}(0)\Psi_{x_i x_i x_i}(\bm r) + \Psi_{x_i}(\bm r)\Psi_{x_i x_i x_i}(0) \notag \\
    & \quad \quad + \Psi_{x_i}(0)\Psi_{x_i x_j x_j}(\bm r) + \Psi_{x_i}(\bm r)\Psi_{x_i x_j x_j}(0) \Big)\Big\rangle.
\end{align}

Here $i = 1,2$ and $j = 1,2$. Then we can compute that 
\begin{align}
    \frac{\sigma^4 \Delta t^2}{16}\Tilde{R_1}(\bm r)&= \frac{\sigma^4 \Delta t^2}{16} \frac{1}{N}\sum_{i=1}^N \frac{1}{|\bm k_i|^2}(k_{i,1}^6 + 3k_{i,1}^4 k_{i,2}^2+ 3k_{i,1}^2 k_{i,2}^4+ k_{i,2}^6 ) \cos(\bm k_i \cdot \bm r)  \notag \\
    &=\frac{\sigma^4 \Delta t^2}{16 N} \sum_{i=1}^N |\bm k_i|^4 \cos(\bm k_i \cdot \bm r). \label{DisSumofR1}
\end{align}
Due to the analysis in \cite{Majda1999}, as $N \to \infty$, this summation \eqref{DisSumofR1} converges, more precisely 
\begin{align}
\lim_{N\to\infty}\frac{\sigma^4 \Delta t^2}{16}\Tilde{R_1}(\bm r)=
   & \frac{\sigma^4 \Delta t^2}{16 }\int_0^L \int_0^{2\pi} k^{2(3-\alpha)} cos(k|\bm r| \cos(\theta)) d\theta d k\\  =& \frac{\sigma^4 \Delta t^2}{16} \int_0^L k^{2(3-\alpha)} J_0(k|\bm r|) d k. \notag
\end{align}
Here $J_0(\cdot)$ is the Bessel function of the first kind.  Similar analysis can be done for $\Tilde{R_2}(\bm r)$ and $\Tilde{R_3}(\bm r)$ and we can obtain that 
\begin{align}
    \lim_{N\to\infty} \frac{\Delta t^2}{4} \Tilde{R_2}(\bm r) &\sim \lim_{N\to\infty} \frac{\Delta t^2}{4 N} \sum_{i=1}^N |\bm k_i|^4 \cos(\bm k_i \cdot \bm r) \notag \\
    & = \frac{\Delta t^2}{4}\int_0^L \int_0^{2\pi} k^{2(3-\alpha)} cos(k|\bm r| \cos(\theta)) d\theta d k\\
    & = \frac{\Delta t^2}{4} \int_0^L k^{2(3-\alpha)} J_0(k|\bm r|) d k. \notag
\end{align}
\begin{align}
    \lim_{N\to\infty} \frac{\sigma^2 \Delta t}{4} \Tilde{R_3}(\bm r) &= \lim_{N\to\infty} \frac{\sigma^2 \Delta t}{2N} \sum_{i=1}^N |\bm k_i|^2 \cos(\bm k_i \cdot \bm r) \notag \\
    &=\frac{\sigma^2 \Delta t}{2}\int_0^{L} \int_0^{2\pi} k^{2(2-\alpha)} cos(k|\bm r| \cos(\theta)) d\theta d k \notag \\
    &= \frac{\sigma^2 \Delta t}{2} \int_0^L k^{2(2-\alpha)} J_0(k|\bm r|) d k. \notag
\end{align}
Then we can observe that for further analysis we need to estimate the moments of the Bessel function of the first kind. Following the analysis in \cite{Wang2016NumericalAF} and \cite{MILNE-THOMSON1945},  we can obtain that 
\begin{align}
    \int_0^L k^{\beta} J_0(k|\bm r|) d k =& \frac{L}{|\bm r|^{\beta}} \big[(\beta-1)J_0(L|\bm r|) s^{(2)}_{\beta-1, -1}(L|\bm r|) + J_1(L|\bm r|) s^{(2)}_{\beta, 0}(L|\bm r|)\big]\notag \\
    &+\frac{2^{\beta} \Gamma(\frac{\beta+1}{2})}{|\bm r|^{\beta + 1}\Gamma(\frac{1-\beta}{2})}.
    \label{EstformomentofJ0}
\end{align}

For $x$ large enough, the Lommel function $s^{(2)}_{n, m}(x)$ has an asymptotic expansion \cite{MILNE-THOMSON1945}, 
\begin{equation}
    s^{(2)}_{n, m}(x) = x^{n-1} \big(1 - \frac{(n-1)^2 + m^2}{x^2} + \frac{\big((n-1)^2 - m^2\big)\big((n-3)^2 - m^2\big)}{x^4} - \cdots\big).
    \label{LommelFct}
\end{equation}

In this case, the Bessel function also has an asymptotic function for $x$ large enough
\begin{equation}
    J_{\alpha} (x) = \sqrt{\frac{2}{\pi x}} (\cos(x - \frac{\alpha \pi}{2} -\frac{\pi}{4} ) + O(\frac{1}{x}) ).
    \label{BesselFct}
\end{equation}

Substitute \eqref{LommelFct} and \eqref{BesselFct} back into \eqref{EstformomentofJ0}, and we can obtain that as $|\bm r| \to \infty$, the leading term in \eqref{EstformomentofJ0} will become
\begin{equation}
\int_0^L k^{\beta} J_0(k|\bm r|) d k\sim
    \left\{
    \begin{aligned}
        &\frac{1}{|\bm r|^{\beta +1}}, 0 < \beta < \frac{1}{2}\\
        &\frac{1}{|\bm r|^{\frac{3}{2}}}, \quad \beta \geq \frac{1}{2}
    \end{aligned}
    \right.
\end{equation}
Then we go back to the analysis of $\Tilde{R_1}(\bm r), \Tilde{R_2}(\bm r)$ and $\Tilde{R_3}(\bm r)$, due to $\frac{1}{4} <\alpha < 1$, so we can obtain that for $|\bm r|$ large enough
\begin{equation}
    \frac{\sigma^4 \Delta t^2}{16}\Tilde{R_1}(\bm r) + \frac{\Delta t^2}{4} \Tilde{R_2}(\bm r) - \frac{\sigma^2\Delta t}{4} \Tilde{R_3}(\bm r) \sim \frac{O(\sigma^2 \Delta t + \Delta t^2)}{|\bm r|^{\frac{3}{2}}}.
\end{equation}

Together with the original correlation function $R_{i,j}(\bm r)$, we can obtain that when $\frac{1}{4} \leq \alpha < 1$, the asymptotic performance of $R^s_{i, j}(\bm r)$ is consistent with the original system, i.e.
\begin{equation}
    R^s_{i, j}(\bm r) \sim \frac{1}{|\bm r|^{2(1-\alpha)}},\ \frac{1}{4} \leq \alpha < 1. 
\end{equation}

Then from Proposition \ref{PropofFann}, we can obtain that for $t$ large enough,
\begin{equation}
    \mathbf{E}[(\bm x_t^S)^2] \sim t^{\frac{2}{2 - \alpha}},\ \frac{1}{4} \leq \alpha < 1. 
\end{equation}
\hfill $\square$
\end{proof}

\section{Numerical experiments}\label{sec:NumericalResults}
\noindent In this section, we present several numerical results computed by the structure-preserving scheme. To start with, we first discuss how to generate random velocity fields $\bm v(\bm x, t)$. Notice that we have suppressed the dependence of the velocity on $\omega$ in $\bm v(\bm x, t)$ for notation simplicity here. In our numerical experiments, we will consider the particles moving in an isotropic velocity field satisfies Assumption \ref{Assforvpro} and \ref{Assforvregular}. 
And we have mentioned that the covariance functions of $\bm v(\bm x, t)$ have the forms like \eqref{CovFct2DS3} and \eqref{CovFct3DS3}. We take the time-correlation function $D(t) = \exp(-\frac{1}{2} \theta_0^2 t^2)$, where $\theta_0 \geq 0$ represents the strength of the decorrelation in time $t$. As for the energy spectrum function $E(k)$, we will consider the spectrums with the following forms:
\begin{align}
    &E_1(k) = \delta(k - k_0),\ \text{2D flows},\label{2DE1}\\
    &E_2(k) = \frac{9}{2} k^3 k_0^{-4} \exp(-\frac{3}{2} k^2/k_0^2),\ \text{2D flows},\label{2DE2}\\ 
    &E_3(k) = \frac{3}{2} \delta(k - k_0),\ \text{3D flows},\label{3DE3}\\
    &E_4(k) = 16 \sqrt{\frac{2}{\pi}} k^4 k_0^{-5} \exp(-2 k^2 / k_0^2),\  \text{3D flows},\label{3DE4}
\end{align}
with a fixed $k_0$. Here $E_1$ and $E_3$ describe excitation confined to a thin shell, in two- and three-
dimensional spaces, respectively, while $E_2$ and $E_4$ describe less concentrated excitations,
in two- and three-dimensional spaces, respectively; see \cite{Kraichnan1970}.

 For the energy spectrum functions $E_1$ and
$E_3$, we choose $\bm k_n$ isotropically distributed on the surface of a sphere and a circle of radius
$k_0$, respectively. For the energy spectrum function $E_2$, each component of $\bm k_n$ is chosen from
a Gaussian distribution of standard deviation $k_0/\sqrt{3}$, while for the energy spectrum function $E_4$, each component of $\bm k_n$ is chosen from a Gaussian distribution of standard deviation $k_0 /2 $. Finally, $\theta_n$ in \eqref{vGene} is chosen from a Gaussian distribution with standard deviation $\theta_0$. In the remaining part of this paper, we will use $D_0$ as a parameter in our numerical experiments, which is equivalent to $\sigma^2 / 2$.

\subsection{Comparison between structure-preserving and Euler-Maruyama}
\noindent In this subsection, we compare the effective diffusivity computed by the structure-preserving scheme and Euler-Maruyama scheme, where the velocity field is generated by \eqref{vGene} with the setting \eqref{uw2Dform}. Here $E(k) = E_1(k)$ defined in \eqref{2DE1}. Theorem 1 of \cite{Fannjiang1996} has explained why $\theta_0 = 0$ and $D_0 = 0$ leads to anomalous diffusion phenomena. The occurrence of the trapping phenomena results in the random velocity field with this setting shall exhibit the sub-diffusion performance.
\begin{figure}[h]
    \centering
    \includegraphics[width = 0.7\textwidth]{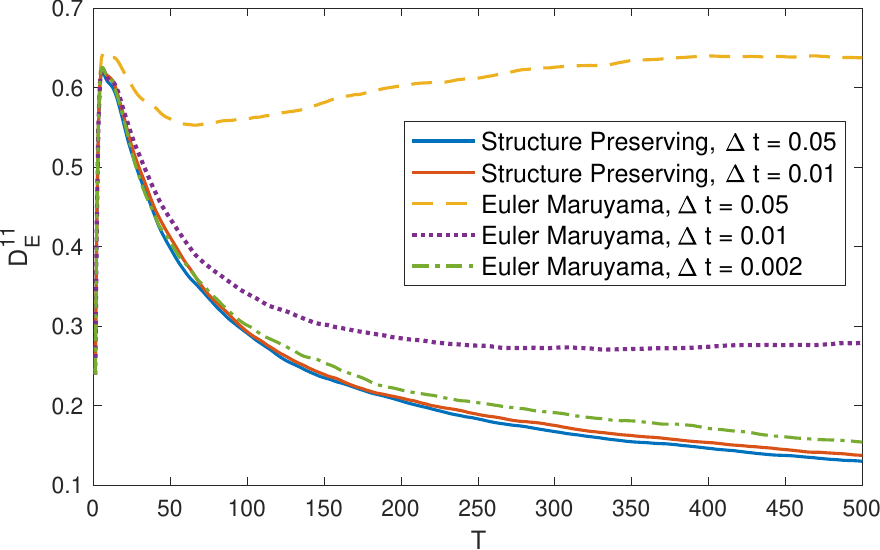}
    \caption{Comparison of $D^{11}_E$ by two numerical schemes with different time step $\Delta t$ in 2D random field with $E_1(k)$ and $\theta_0 = 0, D_0 = 0$}
    \label{2DEulvsSpE1th0}
\end{figure}

\begin{figure}[h]
    \centering
    \includegraphics[width = 0.7\textwidth]{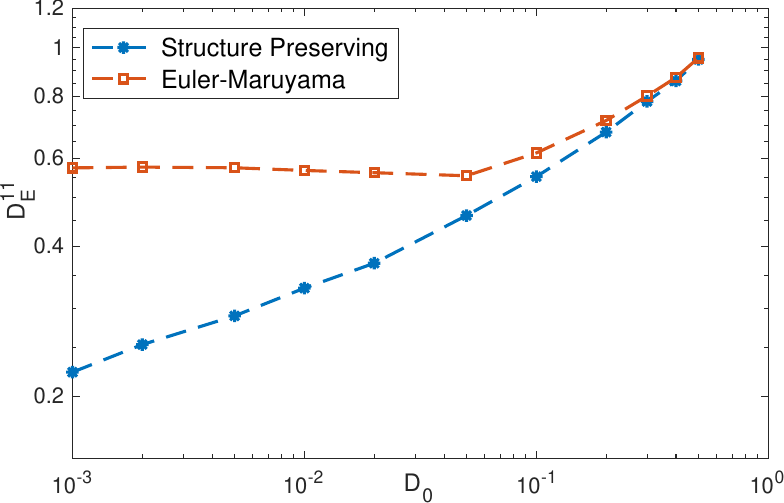}
    \caption{Comparison of $D^{11}_E$  by two numerical schemes with different $D_0$ in 2D random field with $E_1(k)$ and $\theta_0 = 0, \Delta t = 0.1$}
    \label{2DEulvsSpE1th0dt01}
\end{figure}

\begin{figure}[h]
    \centering
    \includegraphics[width = 0.7\textwidth]{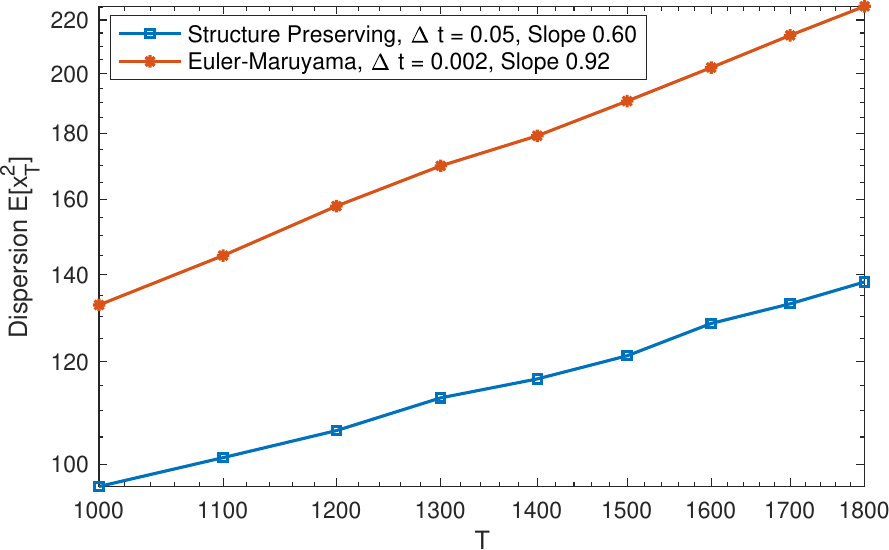}
    \caption{Comparison for dispersion $\mathbb{E}[\bm x_t^2]$, by Log-log plot in 2D random field with $E_1(k)$ and $\theta_0 = 0, D_0 = 0$}
    \label{2DEulvsSpE1Disperth0}
\end{figure}

From Figure \ref{2DEulvsSpE1th0} we can see that, the structure-preserving scheme is not sensitive to the time step $\Delta t$ and their results are very close. they both converge to zero and show the trapping phenomena. While the Euler-Maruyama scheme is sensitive to the time step $\Delta t$ and as we decrease $\Delta t$, it gets closer to the reference result but still converges to a non-zero limit and does not show the sub-diffusion performance. To get a relatively close result, the Euler-Maruyama scheme needs a time step which is much smaller than the structure-preserving scheme. Then we also test for the performance of the Euler-Maruyama scheme under a long time scale and at this time, the Euler-Maruyama scheme will no longer perform sub-diffusion phenomena. This is due to the numerical diffusion phenomena caused by the Euler-Maruyama scheme itself. This is also well illustrated in Figure \ref{2DEulvsSpE1th0dt01}. We can see that when $\Delta t$ is fixed, the effective diffusivity calculated by the Euler-Maruyama method does not decrease with decreasing $D_0$. This phenomenon does not occur in the structure-preserving scheme. This explains that the Euler-Maruyama scheme is very sensitive to the scale relationship between $\Delta t$ and $D_0$. When $D_0 \ll \Delta t $, this scheme drastically changes the behavior of the effective diffusivity. This also supports our Theorem \ref{Coro::subdiffforSP}. To more conveniently tell the diffusion behavior of the random flows, we can observe the dispersion $\mathbb{E}[\bm x_t^2]$ at time $t$. Then Figure \ref{2DEulvsSpE1Disperth0} illustrates that even for a relatively small $\Delta t$, the Euler-Maruyama scheme will approach a diffusion phenomenon in the long time behavior, while the structure-preserving scheme can keep the sub-diffusion phenomenon which is consistent with the theoretical result of the original system.

\subsection{Computing effective diffusivity of 2D and 3D random flows}
\noindent In this subsection, we investigate the convection-enhanced diffusion phenomena and compute the effective diffusivity of two-dimensional incompressible random flows. Here we choose a time step $\Delta t = 0.05$, number of Monte Carlo samples $N_{mc} = 100,000$ and number of Fourier modes $N = 200$ in approximating the random velocity field. For the velocity field generated by \eqref{vGene} according to $E_1(k)$ defined in \eqref{2DE1} and $E_2(k)$ defined in \eqref{2DE2}, we compute their effective diffusivity with different $\theta_0$ and $D_0$ settings.
\begin{figure}[h]
    \centering
    \begin{subfigure}{0.475\textwidth}
    \includegraphics[width = \linewidth]{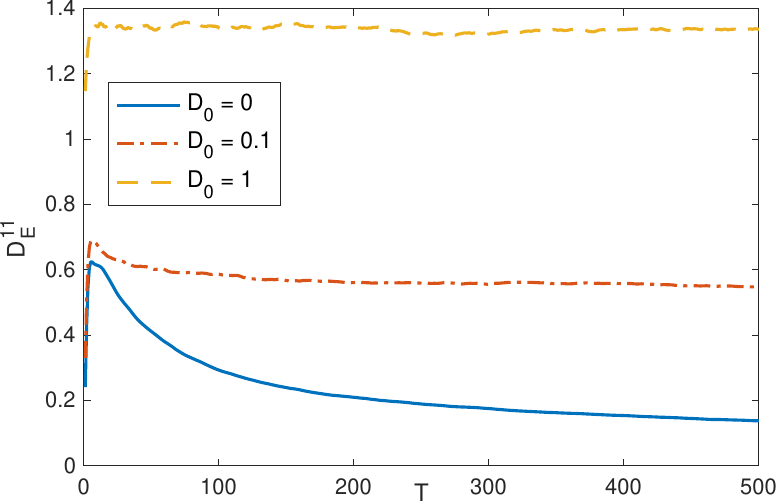}
    \caption{$\theta_0 = 0$}
    \label{2DE1th0}
\end{subfigure}    
\begin{subfigure}{0.48\textwidth}
    \includegraphics[width = \linewidth]{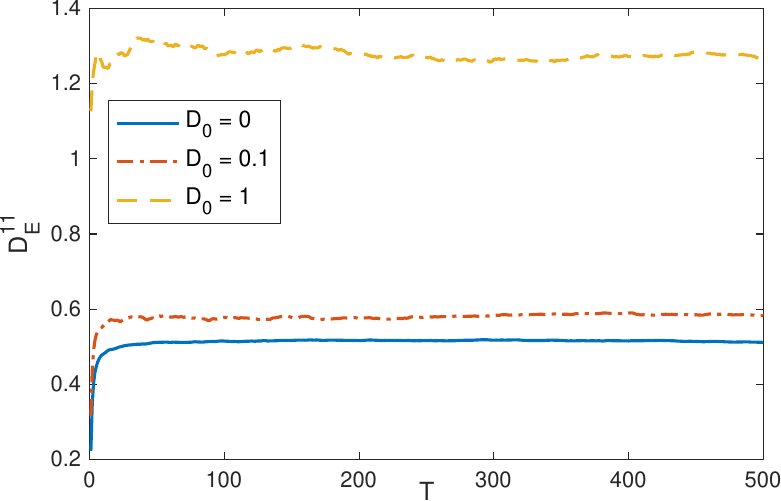}
    \caption{$\theta_0 = 1$}
    \label{2DE1th1}
    
\end{subfigure}    
\caption{Effective diffusivity $D_E^{11}$ v.s. time $t$ in 2D random velocity field with $E_1(k)$}
\label{EffDiff2DE1}
\end{figure}

    

From Figure \ref{EffDiff2DE1} 
, we can see that the effective diffusivity all converges to a non-zero limit except in the case when $\theta_0$ and $D_0$ are both equal to $0$. This implies that for $E_1(k)$ if $D_0 > 0$, its corresponding velocity field will perform normal diffusion phenomena. The energy spectrum density $E_2(k)$ has a similar performance as $E_1(k)$.

\begin{figure}[h]
    \centering
    \begin{subfigure}{0.485\textwidth}
    \includegraphics[width = \linewidth]{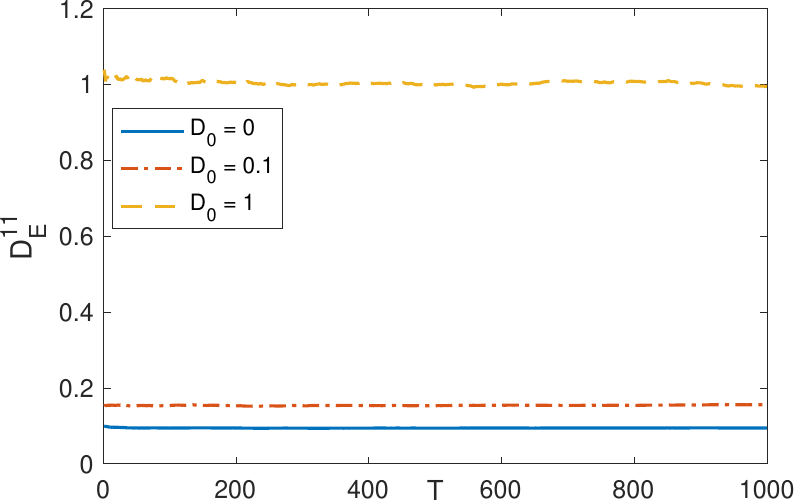}
    \caption{$\theta_0 = 0$}
    \label{3DE3th0}
\end{subfigure}    
\begin{subfigure}{0.47\textwidth}
    \includegraphics[width = \linewidth]{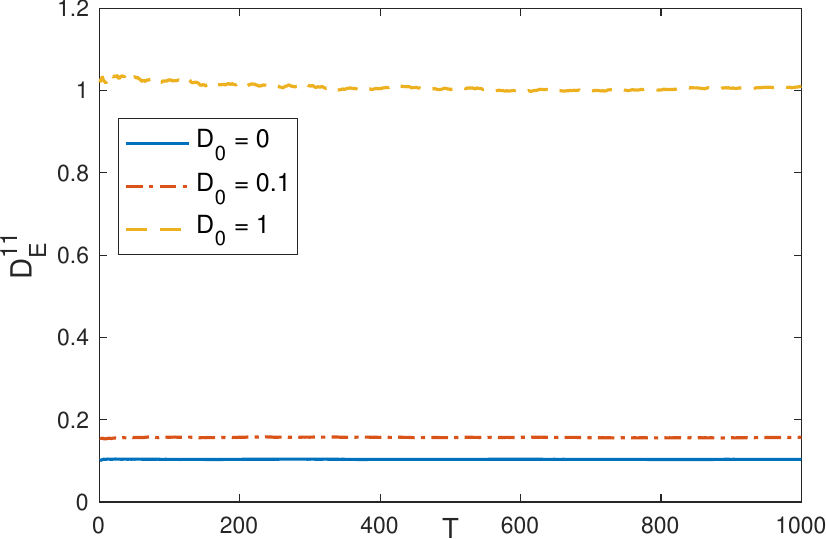}
    \caption{$\theta_0 = 1$}
    \label{3DE3th1}
    
\end{subfigure}    
\caption{Effective diffusivity $D_E^{11}$ v.s. time $t$ in 3D random velocity field with $E_3(k)$}
\label{EffDiff3DE3}
\end{figure}

From Figure \ref{EffDiff3DE3}, we can observe that the effective diffusivity all converges to a non-zero limit in this 3D case. This implies that for $E_3(k)$, its corresponding velocity field will perform normal diffusion phenomena. The energy spectrum $E_4(k)$ performs a similar result to $E_3(k)$. 

\subsection{Anomalous diffusion phenomena in 2D random flows}
\noindent 
In this subsection, we investigate and compute some 2D random velocity fields which will perform anomalous diffusion phenomena. With the condition \eqref{sharpCondi2D}, we can find a $E(k)$ series that does not satisfy this condition.
\begin{align}
&E_5(k) = \sqrt{\frac{6}{\pi}} k_0^{-1} \cdot \exp(-\frac{3k^2}{2k_0^2}), \\
&E_6(k) = \frac{\sqrt[4]{54}}{\Gamma(\frac{3}{4})} k_0^{-1.5} k^{0.5}\cdot \exp(-\frac{3k^2}{2k_0^2}).
\end{align}

Then we plot their dispersion computed by the structure-preserving scheme. Here we set $\theta_0 = 0$ and other parameters remain the same.
\begin{figure}[h]
    \centering
    \includegraphics[width = 0.7\linewidth]{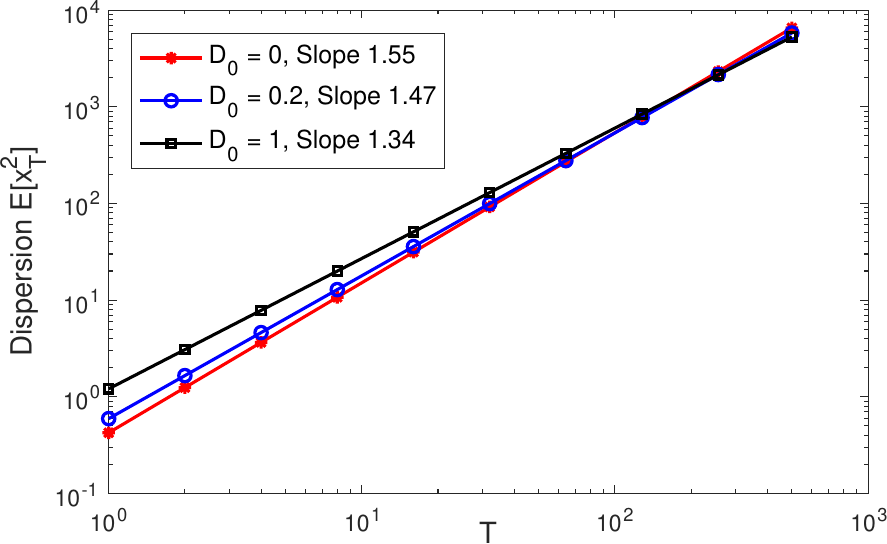}
\caption{Log-log plot of Dispersion $\mathbb{E}[\bm x_t^2]$ v.s. time $t$ in 2D random velocity field with $E_5(k)$}
\label{EffandDis2DE6th0}
\end{figure}

\begin{figure}[h]
    \centering
    \includegraphics[width = 0.7\linewidth]{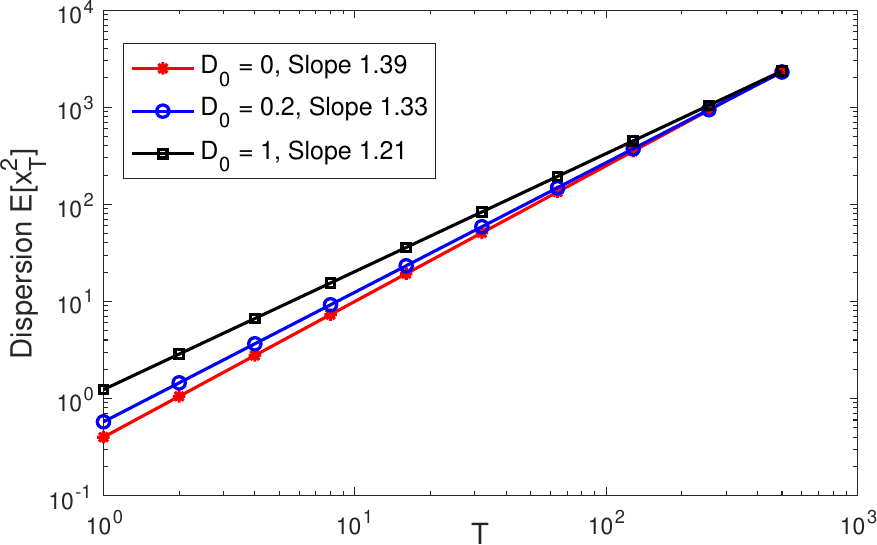}
\caption{Log-log plot of Dispersion $\mathbb{E}[\bm x_t^2]$ v.s. time $t$ in 2D random velocity field with $E_6(k)$}
\label{EffandDis2DE7th0}
\end{figure}

From Figure \ref{EffandDis2DE6th0}, and \ref{EffandDis2DE7th0}, we can see that they both exhibit anomalous diffusion phenomena, specifically super-diffusion phenomena. This is due to the diffusivity is always enhanced in incompressible flows and the computation results are consistent with the theoretical sharp condition \eqref{sharpCondi2D}. Details can be found in \cite{Fannjiang1998}. 

\begin{figure}[h]
    \centering
    \includegraphics[width = 0.7\textwidth]{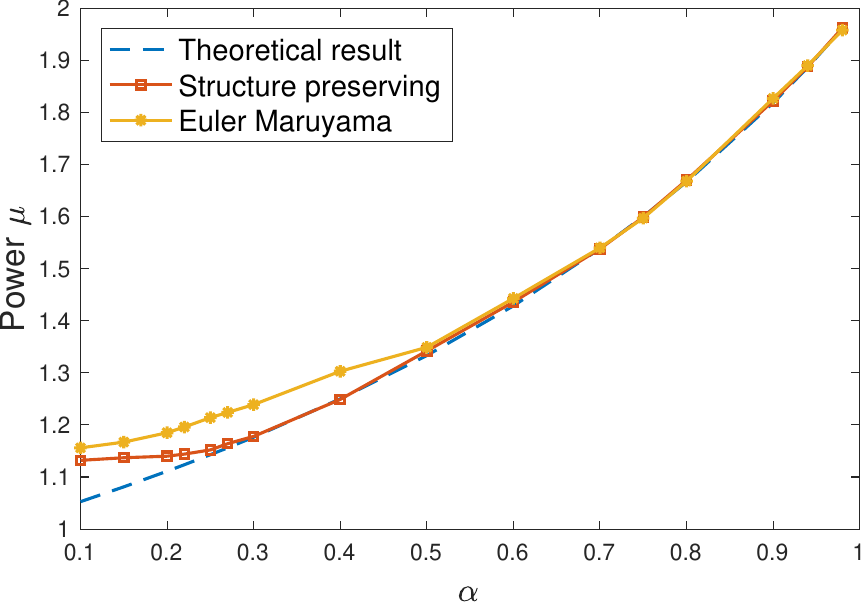}
    \caption{$\mu$ v.s. $\alpha$, comparison for the power law computed by two numerical schemes, $\mathbb{E}[\bm x_t^2] \sim t^{\mu}$, in a 2D random field with spectrum $\hat{R}_{i,j}(\bm k)$ in \eqref{FannVeloSpec}.}
    \label{2DFannEulvsSP}
\end{figure}

 Then, as mentioned in \cite{Fannjiang1998}, we calculate the power law between space and time, $\mu$, of the velocity fields with the correlation functions \eqref{FannVeloSpec} by two different numerical schemes and compare them with the theoretical results. From Figure \ref{2DFannEulvsSP}, we can see that the results of both numerical schemes match the theoretical results well when alpha is relatively large and close to $1$. Meanwhile, when $0.25 < \alpha < 0.5$, the power law computed by the structure-preserving scheme is more in line with the theoretical results. This also coincides with our proof of Theorem \ref{ThmforSuperdiff}.

\subsection{Anomalous diffusion phenomena in 3D random flows}
\noindent 

In this subsection, we investigate and compute some 3D random velocity fields which will perform anomalous diffusion phenomena. With the sharp condition \eqref{sharpCondi3D}, we can find the $E_7(k)$ that does not satisfy this condition as a 3D case.
\begin{equation}
    E_7(k) = 6 k_0^{-2} k \cdot \exp (-\frac{2k^2}{k_0^2}).
    \label{3DE9}
\end{equation}
We can compute its corresponding velocity field's effective diffusivity and dispersion with the same parameters.
\begin{figure}[h]
    \centering
    \includegraphics[width = 0.7\linewidth]{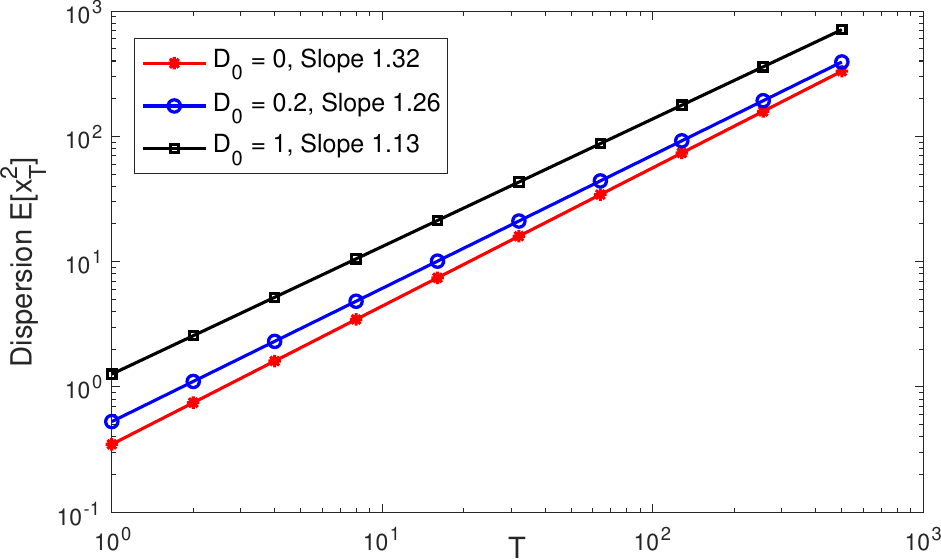}
\caption{Log-log plot of Dispersion $\mathbb{E}[\bm x_t^2]$ v.s. time $t$ in 3D random velocity field with $E_7(k)$}
\label{EffandDis3DE9th0}
\end{figure}

From Figure \ref{EffandDis3DE9th0}, we can see that the 3D random velocity field with energy spectrum $E_7(k)$ performs super-diffusion phenomena which is consistent with the theoretical condition. We also conduct experiments to calculate the power law, $\mu$, of the random field mentioned in \cite{Fannjiang1998} in both numerical schemes in the three-dimensional cases. Now the random velocity field will have the following spectrum 
\begin{equation}
    \hat{R}_{i,j}(\bm k) \sim \frac{1}{|\bm k|^{2\alpha+1}} (\delta_{i,j} - \frac{k_i k_j}{|\bm k|^2}).
        \label{FannVeloSpec3D}
\end{equation}

\begin{figure}[h]
    \centering
    \includegraphics[width = 0.7\textwidth]{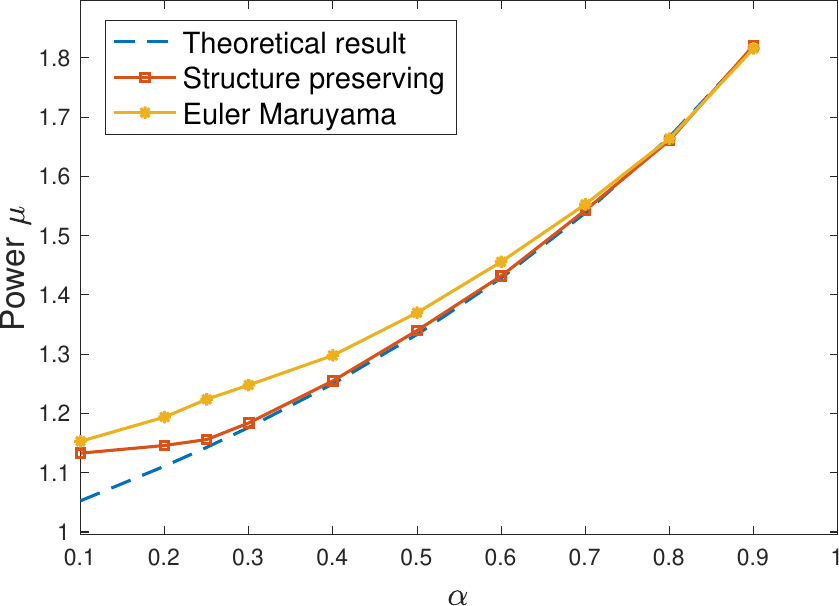}
    \caption{$\mu$ v.s. $\alpha$, comparison for the power law computed by two numerical schemes, $\mathbb{E}[\bm x_t^2] \sim t^{\mu}$, in a 3D random field with spectrum $\hat{R}_{i,j}(\bm k)$ in \eqref{FannVeloSpec3D}.}
    \label{3DFannEulvsSP}
\end{figure}

From Figure \ref{3DFannEulvsSP}, we can see that, when $\alpha$ is not so small, i.e. $\alpha > 0.25$, the power law $\mu$ computed by the structure-preserving scheme matches the theoretical result better than the Euler-Maruyama scheme.

\section{Conclusions}\label{sec:Conclusion}
In this paper, we developed an efficient structure-preserving scheme to compute the effective diffusivity $D^{11}_E$ in 2D and 3D stochastic flows, with convergence analysis based on the Markov semigroup and modified flow system. We also compare the structure-preserving scheme with the Euler-Maruyama scheme and prove that the structure-preserving scheme can keep the original diffusion behavior of the system unchanged in the long time scale while the Euler-Maruyama method cannot. We also conducted numerical experiments to address this point, and the experimental results are consistent with it. For a 2D random field which should exhibit sub-diffusion phenomena when diffusion vanishes, the structure-preserving scheme maintains this property, while the Euler-Maruyama scheme exhibits diffusion phenomena in the long time scale. We also computed the effective diffusivity for some random fields with diffusion phenomena in 2D and 3D. Then based on the criterion in \cite{Fannjiang1996}, we derived and discovered a series of 2D and 3D random fields that would perform anomalous diffusion phenomena. The corresponding numerical experiments also yielded consistent results. Meanwhile, we also prove that for random flows that follow the spectrums \eqref{FannVeloSpec}, \eqref{FannVeloSpec3D} mentioned in \cite{Fannjiang1998}, the structure-preserving scheme can keep the original super-diffusive behavior.
Compared to the Euler-Maruyama scheme  and other Eulerian methods, our structure-preserving particle method has the following advantages: (1) It is scalable in spatial dimension and can handle stochastic 3D problems at ease; (2) It can keep the diffusion behavior of the original system in the long time scale and is not very sensitive to the time step.

\begin{acknowledgements}
\noindent  JX was partially supported by NSF grant DMS-2309520. ZZ was supported by Hong Kong RGC grant (Projects 17300318 and 17307921), National Natural Science Foundation of China  (Project 12171406), Seed Funding Programme for Basic Research (HKU) and Seed Funding for Strategic Interdisciplinary Research Scheme 2021/22 (HKU). 
\end{acknowledgements}

\appendix  

\section{Some formulations and spectral gap conditions in random flows}
\label{AppendixforSomeFormandSGC}
We first define a function space that holds the stationary and ergodic properties in $\mathcal{X}$. Then let $(\mathcal{X}, \mathcal{H}, P_0)$ be a probability space. Here $\mathcal{H}$ is a $\sigma$-algebra over $\mathcal{X}$, and $P_0$ is a probability measure. Then let $\tau_{\bm x}$, $\bm x \in \mathbb{R}^d$ be an almost surely continuous, jointly measurable group of measure-preserving transformation on $\mathcal{X}$ with the following properties: 
\begin{enumerate}
    \item [P1.] $\tau_{\bm 0} = Id_{\mathcal{X}}$, and $\tau_{\bm x + \bm y} = \tau_{\bm x} \tau_{\bm y}, \forall\ \bm x, \bm y \in \mathbb{R}^d$.
    \item [P2.] The mapping $(\Tilde{\bm x}, \bm x) \to \tau_{\bm x} \Tilde{\bm x}$ is jointly measurable.
    \item [P3.] $P_0 (\tau_{\bm x} A) = P_0(A), \forall \bm x\in \mathbb{R}^d, A\in \mathcal{H}$.
    \item [P4.] $\lim_{\bm x \to \bm 0} P_0 (\Tilde{\bm x}: |f \circ \tau_{\bm x}(\Tilde{\bm x}) - f(\Tilde{\bm x})| \geq \epsilon ) = 0,\ \forall f\in L^2(\mathcal{X}), \epsilon > 0$.
    \item [P5.] If $P_0 (A \Delta \tau_{\bm x}A ) = 0, \forall \bm x \in \mathbb{R}^d$, then $A$ is a trivial event.
\end{enumerate}

Then $\tau_{\bm x}$ induces a strongly continuous group of unitary mapping $U^{\bm x}$ on $L^2(\mathcal{X})$, which satisfies
\begin{equation}
    U^{\bm x}f(\Tilde{\bm x}) = f(\tau_{\bm x}(\Tilde{\bm x})),\ f\in L^2(\mathcal{X}),\ \bm x\in \mathbb{R}^d.
\end{equation}

Next, we include the time variable $t$ and investigate the Markovian property. Let $\Omega$ be the space of $\mathcal{X}$-valued continuous function $C([0,\infty); \mathcal{X})$ and let $\mathcal{B}$ be its Borel $\sigma$-algebra. Let $P^t, t\geq 0$, be a strongly continuous Markov semigroup on $L^2(\mathcal{X})$, which satisfies the following properties.
\begin{enumerate}
    \item[Q1.] $P^t \textbf{1} = \textbf{1}$ and $P^t f \geq 0$ if $f\geq 0$.
    \item[Q2.] $\int P^t f dP_0 = \int f dP_0$, for all $f \in L^2(\mathcal{X}), t \geq 0$.
    \item[Q3.] $\textbf{E}_{\Tilde{\bm x}}[f(\theta_{t+h}(\omega))| \mathcal{B}_{\leq t}] = P^h F(\omega(t))$, where $F(\Tilde{\bm x}):= \textbf{E}_{\Tilde{\bm x}}f$ for any $f\in L^1(\Omega)$, $ t,h\geq 0, \Tilde{\bm x} \in \mathcal{X}$.
\end{enumerate}

Here $\mathcal{B}_{\leq t}$ are the $\sigma$-algebras generated by measurable events preceding time $t$, and $\theta_t(\omega(\cdot)) = \omega(\cdot + t), t\geq 0$ is an standard shift operator on the path space $(\Omega, \mathcal{B})$.
Then a stationary measure $P$ can be defined on the path space $(\Omega, \mathcal{B})$ by
\begin{equation}
    P(B) = \int P_{\Tilde{\bm x}}(B) P_0(d\Tilde{\bm x}) , B \in \mathcal{B}, \notag
\end{equation}
and we use $\textbf{E},\ \textbf{E}_{\Tilde{\bm x}}$ to denote the corresponding expectation operators taken over the measure $P,\  P_{\Tilde{\bm x}}$. Here $P_{\Tilde{\bm x}}$ can be considered as the conditional probability for all events in $\mathcal{B}$ condition with the initial point lies on $\Tilde{\bm x}$. We denote $\theta_t$ as a standard shift operator on the path space $(\Omega, \mathcal{B})$.   
\begin{proposition}
    $P$ is invariant under the action of $\theta_t$ and $\tau_{\bm x}$ for any $(t, \bm x) \in \mathbb{R}^+ \times \mathbb{R}^d$.
\end{proposition}

Let $L: \mathcal{D}(L) \to L^2(\mathcal{X})$ be the generator of the semigroup $P^t$. We assume that the generator $L$ satisfies the following spectral gap condition 
\begin{equation}
    -(Lf, f)_{L^2(\mathcal{X})} \geq c_1 ||f||^2_{L^2(\mathcal{X})},\ c_1 > 0.
    \label{SpecGapCondiL}
\end{equation}

The spectral gap condition \eqref{SpecGapCondiL} guarantees the exponential decay property
\begin{equation}
    ||P^t f||_{L^2(\mathcal{X})}\leq \exp(-2c_1 t)||f||_{L^2(\mathcal{X})},\ c_1 > 0,\ f\in L_0^2(\mathcal{X}).
    \label{SpectralGapProp}
\end{equation}

See \cite{Doukhan1994}, \cite{Rosen2012}. The exponential decay property \eqref{SpectralGapProp} for $P^t$ is equivalent to $\varrho$-mixing of the process $X(t), t \geq 0$ and contributes significantly to proving the existence of the effective diffusivity.

\section{The behavior of the stream function $\Psi$ in two numerical schemes}
\label{AppendixforStreamfct}

\begin{lemma}
\label{lemma::streamfctcon}
    Let $\Psi$ defined in \eqref{Hamiltion} be the stream function of a 2D divergence-free random velocity field. 
     We can treat the mean value of the Hamiltonian $\Psi$ as a function of time $t$. Then it can be obtained that
    \begin{equation}
    \frac{d \Psi}{dt} =  \frac{\sigma^2 }{2} \Delta\Psi + \sigma \nabla \Psi\cdot \frac{d \bm w_t}{dt}.
    \label{EquofLemm3}
\end{equation}
    \label{Lemmaforphicontin}
\end{lemma}

\begin{proof}
    We apply Ito’s formula to the stream function $\Psi$ for $\bm x$ satisfies the original SDEs \eqref{PasTraSDEforVx}, to find 
\begin{equation}
    d\Psi = [(\nabla \Psi)\cdot \bm v + \frac{1}{2} \text{Tr}((\sigma \mathbf{I}_2) (\mathbf{H}_{\Psi}) (\sigma \mathbf{I}_2))] dt + (\nabla \Psi)\cdot (\sigma \mathbf{I}_2) d\bm w_t,
\end{equation}
here $\mathbf{H}_{\Psi}$ is the Hessian matrix of $\Psi$ with respect to $\bm x$. Then due to $\bm v = \nabla^{\perp} \Psi$ is divergence-free, we can compute that 
\begin{equation}
    (\nabla \Psi)\cdot \bm v = (\nabla \Psi)\cdot (\nabla^{\perp} \Psi) =0,\quad \frac{1}{2} \text{Tr}((\sigma \mathbf{I}_2) (\mathbf{H}_{\Psi}) (\sigma \mathbf{I}_2)) = \frac{\sigma^2}{ 2} \Delta \Psi.
\end{equation}

Then we can obtain that
\begin{equation}
    \frac{d \Psi}{dt} = \frac{\sigma^2}{2} \Delta \Psi + \sigma \nabla \Psi\cdot \frac{d \bm w_t}{dt}.
    \label{PsiafterIto}
\end{equation}
The integral of $\sigma \nabla \Psi\cdot \frac{d \bm w_t}{dt}$ is a mean zero martingale, i.e. $\mathbf{E} [\int_0^t \sigma \nabla \Psi\cdot d \bm w_t] = 0$. 

\hfill $\square$
\end{proof}

\begin{lemma}
\label{lemma::streammodified}
    Let $\Psi^E$ and $\Psi^S$ be the stream function of a 2D divergence-free random velocity field satisfying the modified equation systems of Euler-Maruyama scheme and structure-preserving scheme correspondingly. 
    The Euler-Maruyama method will introduce an extra numerical diffusion term in the behavior of the mean value of $\Psi^E$, i.e.
    \begin{align}
    \frac{d\Psi^E}{d t} &= \frac{\Delta t}{2  N \sqrt{N}} \sum_{n=1}^N M_n \Delta \Psi_n +\frac{\sigma^2 }{2}  M_{\sigma} \Delta \Psi^E + \nabla \Psi^E \cdot \mathbf{G}_E \frac{d \bm w_t}{dt}.
    \label{EquofLem44}
    \end{align} 
    Here 
    \begin{align}
        M_n &= \frac{\big( \sum_{i=1}^N (k_{n,1} k_{i,2} - k_{n,2} k_{i,1})(\xi_i \cos(\bm k_i \cdot \bm x) + \zeta_i \sin(\bm k_i \cdot \bm x))\big)^2}{ k_{n,1}^2 + k_{n,2}^2}\geq 0, \notag \\
        M_{\sigma} &= 1 + \frac{3 \Delta t^2}{4} \Psi_{x_1 x_2}^2 + \frac{\Delta t^2}{4}(\Psi_{x_1 x_1}^2 - \Psi_{x_1 x_1}\Psi_{x_2 x_2}+ \Psi_{x_2 x_2}^2 ) \geq 0, \notag
    \end{align}
    and the integral of $\nabla \Psi^E \cdot \mathbf{G}_E \frac{d \bm w_t}{dt}$ is a mean zero martingale, i.e.  $\mathbb{E} [\int_0^t \nabla \Psi^E \cdot \mathbf{G}_E d \bm w_t]$ $ = 0$. Compared with \eqref{EquofLemm3} in Lemma \ref{lemma::streamfctcon}, the Euler-Maruyama scheme will introduce an additional numerical diffusion term. While the structure-preserving scheme will not introduce this error term, i.e.
    \begin{align}
    \frac{d\Psi^S}{d t} &= \frac{\sigma^2 }{2} M_{\sigma} \Delta \Psi^S + \nabla \Psi^S \cdot \mathbf{G}_S \frac{d \bm w_t}{dt}.
    \end{align}
    \label{Lemma4.1}
\end{lemma}

\begin{proof}
Inspired by the idea in \cite{PAVLIOTIS20091030}, we apply Ito's formula to the stream function $\Psi^E (\bm x)$ for the random flow for $\bm x$ satisfies the modified equation \eqref{ModifiedEquofEul} of the Euler-Maruyama method. We can find that 
\begin{equation}
    d \Psi^E = [(\nabla \Psi^E)\cdot (\bm v+\Delta t \Tilde{\bm v}_E) + \frac{\sigma^2}{2} \text{Tr}(\mathbf{G}_E^T (\mathbf{H}_{\Psi}) \mathbf{G}_E)] dt + (\nabla \Psi^E)\cdot \mathbf{G}_E d\bm w_t. 
\end{equation}

Then we can compute that,
\begin{align}
    (\nabla \Psi^E)\cdot (&\bm v+\Delta t \Tilde{\bm v}_E) = \frac{\Delta t}{2} (\Psi_{x_1}^2 \Psi_{x_2 x_2} + \Psi_{x_2}^2 \Psi_{x_1 x_1} - 2 \Psi_{x_1} \Psi_{x_2} \Psi_{x_1 x_2}) \notag \\
    &\quad + \frac{\sigma^2 \Delta t}{4} (\Psi_{x_1} \Psi_{x_1 x_1 x_2} + \Psi_{x_1} \Psi_{x_2 x_2 x_2} - \Psi_{x_2} \Psi_{x_1 x_2 x_2} - \Psi_{x_2} \Psi_{x_1 x_1 x_1} ).
    \label{ModifiedEQterm}
\end{align}
\begin{align}
    \frac{\sigma^2}{2} \text{Tr}(\mathbf{G}_E^T &(\mathbf{H}_{\Psi}) \mathbf{G}_E)  = \frac{\sigma^2}{2} M_2  \Delta \Psi^E = -\frac{\sigma^2 k_0^2}{2} M_2 \Psi^E.
    \label{eq431}
\end{align}
Here
\begin{equation}
    M_2 = 1 + \frac{3 \Delta t^2}{4} \Psi_{x_1 x_2}^2 + \frac{\Delta t^2}{4}(\Psi_{x_1 x_1}^2 - \Psi_{x_1 x_1}\Psi_{x_2 x_2}+ \Psi_{x_2 x_2}^2) \geq 1. \notag
\end{equation}

And for the two terms in \eqref{ModifiedEQterm} for the case that the flow with the stream function \eqref{Hamiltion}, we can obtain that 
\begin{align}
    &\frac{\Delta t}{2} (\Psi_{x_1}^2 \Psi_{x_2 x_2} + \Psi_{x_2}^2 \Psi_{x_1 x_1} - 2 \Psi_{x_1} \Psi_{x_2} \Psi_{x_1 x_2}) = \notag \\ 
    &\frac{-\Delta t}{2 N \sqrt{N}} \sum_{n=1}^N \big( \sum_{i=1}^N (k_{n,1} k_{i,2} - k_{n,2} k_{i,1})(\xi_i \cos(\bm k_i \cdot \bm x) + \zeta_i \sin(\bm k_i \cdot \bm x))\big)^2 \Psi_n.
    \label{eq432}
\end{align}

And we can easily verify that
\begin{equation}
    \frac{\sigma^2 \Delta t}{4} (\Psi_{x_1} \Psi_{x_1 x_1 x_2} + \Psi_{x_1} \Psi_{x_2 x_2 x_2} - \Psi_{x_2} \Psi_{x_1 x_2 x_2} - \Psi_{x_2} \Psi_{x_1 x_1 x_1} ) = 0.
    \label{eq433}
\end{equation}

Then together with \eqref{eq431}, \eqref{eq432} and \eqref{eq433}, we can obtain that
\begin{align}
    \frac{d \Psi^E}{d t} &= \frac{\Delta t}{2 N \sqrt{N}} \sum_{n=1}^N M_{n} \Delta \Psi_n + \frac{\sigma^2 }{2} M_2 \Delta \Psi^E + \nabla \Psi^E\cdot \mathbf{G}_E \frac{d \bm w_t}{dt}  \\
    \label{dStreamFctforEulscheme}
\end{align}
where  $\mathbf{E} [\int_0^t \nabla \Psi^E\cdot \mathbf{G}_E d \bm w_t] = 0$ and
\begin{equation}
    M_n = \big(\sum_{i=1}^N (k_{n,1} k_{i,2} - k_{n,2} k_{i,1})(\xi_i \cos(\bm k_i \cdot \bm x) + \zeta_i \sin(\bm k_i \cdot \bm x))\big)^2 \Big / (k_{n,1}^2 + k_{n,2}^2) \geq 0
\end{equation}
A similar analysis can be done for the structure-preserving scheme, and then we can finally obtain that 
\begin{align}
    \frac{d \Psi^S}{d t} &= \frac{\sigma^2 }{2} M_2 \Delta \Psi^S + \nabla \Psi^S\cdot \mathbf{G}_S \frac{d \bm w_t}{dt} 
    \label{dStreamFctforSypscheme}
\end{align}

\hfill $\square$
\end{proof}

\begin{remark}
    The first term on the right-hand side of \eqref{EquofLem44} does not depend on the diffusion coefficient $\sigma$ but on the time step $\Delta t$. When $\sigma^2 \ll\Delta t$, it drastically changes the behavior of  $\mathbf{E}[\Psi(t)]$ compared with \eqref{EquofLemm3}.
\end{remark}
\begin{remark}
    For energy spectral densities that are distributed on the same unit circle as $E_1(k) = \delta(k-k_0)$, we have $k_{n,1}^2 + k_{n,2}^2 = k_0^2$. Then we can further obtain that 
    \begin{equation}
    \frac{d \Psi}{dt} = - \frac{\sigma^2 k_0^2}{2} \Psi + \sigma \nabla \Psi\cdot \frac{d \bm w_t}{dt},
    \end{equation}
    which means the mean value of the stream function, $\mathbf{E}[\Psi(t)]$ will decay like $e^{-2t/(\sigma^2 k_0^2)}$, i.e. $\mathbf{E}[\Psi(t)] \approx  e^{-2t/(\sigma^2 k_0^2)}\mathbf{E}[\Psi(0)]$.
\end{remark}

%
%

\bibliographystyle{spmpsci}      
\bibliography{Reference}   

%
%

\end{document}